\documentclass{article}
\usepackage[english]{babel}
\usepackage[utf8]{inputenc}
\usepackage[style=numeric,maxbibnames=99]{biblatex}
\usepackage{amssymb}
\usepackage{amsmath}
\usepackage{graphicx,subcaption,float,epstopdf}
\usepackage{amsthm}
\usepackage{xcolor}
\usepackage{hyperref}
\usepackage{authblk}
\usepackage{comment}
\numberwithin{equation}{section}

\usepackage{mathtools}

\newtheorem{theorem}{Theorem}[section]
\newtheorem{proposition}[theorem]{Proposition}

\newtheorem{definition}[theorem]{Definition}

\newtheorem{example}[theorem]{Example}

\usepackage{tikz}
\usetikzlibrary{decorations.pathreplacing}
\usetikzlibrary{positioning}
\usetikzlibrary{arrows}
\usetikzlibrary{arrows.meta}
\usetikzlibrary{calc}
\usetikzlibrary{shapes}

\graphicspath{ {./figures/} }

\usepackage{geometry}
\newlength\titlebox 
\newcommand{\mattia}[1]{\textcolor{blue}{#1}}
\newcommand{\olga}[1]{\textcolor{red}{#1}}
\newcommand{\jim}[1]{\textcolor{violet}{#1}}

\title{Recurrent sets and Lyapunov functions: \\
from topological dynamics to weak KAM theory}

\author[1]{Olga Bernardi}
\author[2]{Mattia Sensi}
\author[3]{Jim Wiseman}
\affil[1]{Department of Mathematics, University of Padova, Italy}
\affil[2]{Department of Mathematics, University of Trento, Italy}
\affil[3]{Department of Mathematics, Agnes Scott College, Decatur, GA, USA}

\date{\today}

\begin{document}
\maketitle

\begin{abstract} 
\noindent The aim of this survey is to clarify the intimate relations between Lyapunov functions and chain recurrent sets. The study of this subject comes from a seminal paper by Conley and has had recent advances by the so-called “Weak KAM Theory for homeomorphisms”. We provide a survey of the state of the art, presenting recent improvements and connections with the stationary Hamilton--Jacobi equation. Finally, in order to illustrate the distinctions between recurrent sets in a non-trivial case as well as their dependence on a parameter, we discuss a compartmental model from mathematical epidemiology and a fast-slow version of a prey-predator model with Allee effect.
\end{abstract}

\section{Introduction}

One of the basic problems in the qualitative study of a dynamical system is to distinguish recurrent motions from those trajectories which move irreversibly towards simpler asymptotic states. Lyapunov functions provide a natural tool for this purpose: they are non-increasing along trajectories and therefore identify the regions where the dynamics has a preferred direction. Recurrence, on the other hand, describes the possibility for a trajectory to return, exactly or approximately, to a previously visited state. Understanding the relation between these two viewpoints is the main theme of this survey.

A fundamental contribution to this problem is due to C.C. Conley. His Fundamental Theorem of Dynamical Systems shows that every continuous flow admits a complete Lyapunov function \cite{Conley1978,conley1972gradient}. This function is constant on the chain recurrent part of the dynamics and strictly decreases along all other trajectories. It therefore decomposes the phase space into two qualitatively different regions: a recurrent region and a gradient-like region. Equivalently, the neutral set of the Lyapunov function, namely the set where no strict decrease occurs, coincides with the chain recurrent set. This point of view gives chain recurrence a distinguished role in topological dynamics and provides a global decomposition which does not depend on the particular metric chosen on the compact phase space.

There are many other notions of recurrence, however, and further work showed  that the chain recurrent set is not always the smallest dynamically meaningful set outside which a Lyapunov function may decrease. A related line of research had already been initiated by J. Auslander through the notion of generalized recurrence, defined in terms of all continuous Lyapunov functions
\cite{auslander1964generalized}. Later, R. Easton introduced strong chain recurrence by strengthening the usual chain condition \cite{easton1978chain}. Instead of requiring each individual jump in a pseudo-orbit to be small, one requires the sum of all jumps to be arbitrarily small. The resulting strong chain recurrent set may be strictly smaller than the chain recurrent set and, unlike the latter, may depend on the metric. These notions therefore provide increasingly refined ways of separating recurrent behaviour from gradient-like dynamics.

A new perspective on these questions emerged from weak KAM theory, where KAM stands for Kolmogorov--Arnold--Moser. A. Fathi and P. Pageault showed that variational techniques originating in Aubry--Mather theory can be used to construct Lipschitz-continuous Lyapunov functions whose neutral set coincides with the strong chain recurrent set for a homeomorphism \cite{Fathi_Pageault_2015_AubryMather,pageault2009conley}. Their construction associates a cost with chains and turns the problem of recurrence into a variational minimization problem. This point of view was subsequently extended to continuous flows, leading to Conley-type Lyapunov theorems for strong chain recurrence \cite{bernardi2019conley,bernardi2019existence,bernardi2021conley}.
Thus, ideas developed originally in Lagrangian and Hamiltonian dynamics entered naturally into a problem arising in topological dynamics.

The connection with Hamiltonian dynamics becomes even more transparent when one considers the stationary Hamilton--Jacobi equation. Given a smooth vector field \(V\) on a compact Riemannian manifold, one can associate with it the Ma\~n\'e Hamiltonian
\[
H_V(x,p)=\frac12\|p\|_x^2+p(V(x)).
\]
The corresponding projected Aubry set \(\mathcal A_V\) provides an intermediate object between strong chain recurrence and chain recurrence: under the hypotheses considered in this paper,
\[
\mathcal{SCR}(\phi_t)
\subseteq
\mathcal A_V
\subseteq
\mathcal{CR}(\phi_t).
\]
Moreover, critical subsolutions of the Hamilton--Jacobi equation give rise to Lyapunov functions for the original flow. In particular, weak KAM theory provides \(C^{1,1}\) critical subsolutions which are strict outside the projected Aubry set
\cite{Bernard2007,Bernard2010,Fathi2012}. This makes it possible to compare, within a single framework, objects coming from topological dynamics, Lyapunov theory, and Hamilton--Jacobi equations.

The present survey is intended as a guided and accessible introduction to Conley's theory of recurrence and Lyapunov functions, with particular attention to the existence and regularity of Lyapunov functions and to the more recent developments inspired by weak KAM theory. It is intended for a broad audience. For this reason, we emphasize the geometric and dynamical meaning of the main notions, rather than their most general abstract formulation, and we accompany the theoretical results with elementary examples and applications.

Our purpose is not to provide an exhaustive treatment of any of the theories involved. Rather, we aim to explain the common mechanisms behind the different notions of recurrence and the corresponding Lyapunov functions. Particular attention is devoted to the passage from Conley's Lyapunov function to strong chain recurrence, to the variational constructions introduced by Fathi and Pageault, and to the relation between neutral sets, generalized recurrence, and weak KAM theory.

The paper is organized as follows. In Section~\ref{Charles Conley}, we recall chain recurrence and Conley's Fundamental Theorem of Dynamical Systems, including the attractor-repeller decomposition underlying the construction of a complete Lyapunov function. In Section~\ref{sec:strength}, we introduce strong chain recurrence and review Conley-type Lyapunov theorems for homeomorphisms and flows, emphasizing the variational approach inspired by weak KAM theory. Section~\ref{sec:neut_sets} discusses generalized recurrence and its relation with the neutral sets of continuous and Lipschitz-continuous Lyapunov functions. In Section~\ref{sec:ham_jac}, we turn to the stationary Hamilton--Jacobi equation, introducing the Ma\~n\'e Hamiltonian and the projected Aubry set and explaining their relation with chain and strong chain recurrence.

Finally, Section~\ref{sec:examples} illustrates these ideas through models from mathematical biology. We first consider a reduced SIRWS
(Susceptible--Infected--Recovered--Waning immunity--Susceptible) compartmental epidemic model \cite{dafilis2014dynamical,jardon2021geometric}. In a singular regime, this model provides a concrete mechanism for the strict inclusion
\[
\mathcal{SCR}(\phi_t)
\subsetneq
\mathcal{CR}(\phi_t).
\]
We then study a planar fast-slow predator-prey model with Allee effect \cite{borsotti2026geometric} and describe how its recurrent sets change as a parameter varies, passing through regimes characterized by equilibria, a heteroclinic cycle, and a periodic orbit. These examples are intended to show that the notions discussed in the preceding sections are not merely abstract refinements, but can provide useful information in the qualitative analysis of concrete nonlinear systems.

\section{The Fundamental Theorem of Dynamical Systems} \label{Charles Conley}
\begin{figure}[h!]
\centering
  \centering
\begin{tikzpicture}
 \node at (0,0) {\includegraphics[width=.35\linewidth]{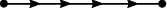}};
\node at (-2.7,0.4) {$0$};
\node at (2.7,0.4) {$\pi$};
  \end{tikzpicture}
\caption{The phase portrait of Example \ref{primissimo}.}
\label{diagramma fase}
\end{figure}

Throughout the exposition, $(X,d)$ is a compact metric space, $f:X \to X$ is a homeomorphism and $\phi_t: X \times \mathbb{R} \to X$ is a continuous flow. We begin with a simple motivating example.
\begin{example} \label{primissimo}
On the interval $[0,\pi]$, let $\phi_t$ be the flow generated by the vector field $V(x)=\sin x$, which is the negative gradient of the function $u(x)=\cos x$, since $V(x)=-u'(x)$. The only recurrent points are the fixed points $0$ and $\pi$. Every other point lies on a heteroclinic orbit connecting $0$ to $\pi$, along which $u$ is strictly decreasing. Indeed,
\[
\frac{d}{dt}u(\phi_t(x))
=u'(\phi_t(x))\,V(\phi_t(x))
=-\sin^2(\phi_t(x))<0
\]
for $x\in(0,\pi)$. Hence
\[
u(\phi_t(x))<u(x),\qquad t>0.
\]
We refer to Figure \ref{diagramma fase} for a visualization.
\end{example}

\noindent Conley's Fundamental Theorem of Dynamical Systems shows that for every system, the dynamics decomposes in a similar way, into a recurrent part and a gradient-like part connecting the recurrent parts. \\
We recall some basic notions in order to state the theorem. First, the notions of Lyapunov function and neutral set, both in the discrete and in the continuous case. 

\begin{definition} $(i)$ Let $f: X \to X$ be a homeomorphism. $u \in C^0(X,\mathbb{R})$ is a \emph{Lyapunov function} for $f$ if $u(f(x)) \le u(x) \quad \forall x \in X$. In such a case, the corresponding neutral set is
$$N(u) = \{ x \in X: \ u(f(x)) = u(x)\}.$$
In particular, $u$ is a \emph{first integral} for $f$ if and only if $N(u) = X$. \\
$(ii)$ Let $\phi_t: X \times \mathbb{R} \to X$ be a continuous flow. $u \in C^0(X,\mathbb{R})$ is a \emph{Lyapunov function} for $\phi_t$ if $u(\phi_t(x)) \le u(x) \quad \forall x \in X \textrm{ and } t \ge 0$. In such a case, the corresponding neutral set is
$$N(u) = \{ x \in X: \ \exists t > 0 \textrm{ such that } u(\phi_t(x)) = u(x)\}.$$
In particular, $u$ is a \emph{first integral} for $\phi_t$ if and only if $N(u) = X$.
\end{definition}
\noindent We remark that ``first integral'' is often replaced with ``constant of motion''. We now recall the definition of chain recurrent set. 
\begin{definition} \label{CR DEF}
$(i)$ Let $f: X \to X$ be a homeomorphism. A point $x \in X$ is said to be \emph{chain recurrent} for $f$ if for any $\varepsilon > 0$ there exists a finite sequence $(x_i)_{i=1,\dots,n} \subset X$ such that $x_1 =x$ and setting $x_{n+1} :=x$, we have 
\begin{equation} \label{omeo}
d(f(x_i),x_{i+1}) < \varepsilon, \qquad i= 1,\dots,n. 
\end{equation}
$(ii)$ Let $\phi_t: X \times \mathbb{R} \to X$ be a continuous flow. A point $x \in X$ is said to be \emph{chain recurrent} for $\phi_t$ if for any $\varepsilon > 0$ and $T > 0$, there exists a finite sequence $(x_i,t_i)_{i=1,\dots,n} \subset X \times \mathbb{R}$ such that $t_i \ge T$ for all $i$, $x_1 =x$ and setting $x_{n+1} :=x$, we have 
\begin{equation} \label{flusso}
d(\phi_{t_i}(x_i),x_{i+1}) < \varepsilon, \qquad i= 1,\dots,n.
\end{equation}
\end{definition}

\begin{figure}[h!]
\centering
  \centering
\begin{tikzpicture}
 \node at (0,0) {\includegraphics[width=.35\linewidth]{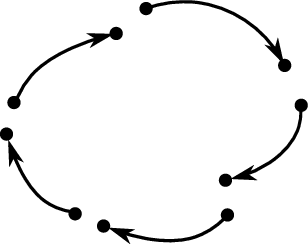}};
\node at (-0.3,2.4) {$x_1$};
\node at (3,1.2) {$f(x_1)$};
\node at (3.2,0.3) {$x_2$};
\node at (1,-0.7) {$f(x_2)$};
\node at (1.7,-1.6) {$x_3$};
\node at (-0.7,-1.5) {$f(x_3)$};
\node at (-1.5,-2) {$x_4$};
\node at (-3.3,-0.2) {$f(x_4)$};
\node at (-2.9,0.4) {$x_5$};
\node at (-0.8,2) {$f(x_5)$};
  \end{tikzpicture}
\caption{A pseudo-orbit for $f$ from $x_1$ to itself.
\label{fig:fig1}}\end{figure}

\noindent The finite sequences $(x_i)_{i=1,\dots,n}$ and $(x_i,t_i)_{i=1,\dots,n}$ in the previous definition are called \emph{pseudo-orbits}: they essentially resemble a true orbit when positions are known only to within a finite accuracy $\varepsilon$. We refer to Figure \ref{fig:fig1} for a visualization. In what follows, the set of chain recurrent points is denoted by $\mathcal{CR}$. The chain recurrent set is easily seen to be closed, invariant, and independent of the metric $d$. For compact ambient spaces, $\mathcal{CR}$ is in particular compact and non-empty. We refer e.g. to \cite[pp.
12–13, 71 and 109–110]{akin1993general}, \cite[Thm. 3.3B]{conley1972gradient} and \cite[p. 248]{zheng2000attractor} for all details. 
\\ \\
\noindent The study of the intimate relation between the chain recurrent set and Lyapunov functions comes from the celebrated paper \cite{Conley1978} by C.C. Conley where it is proved -- see the statement  below -- that parts of the space are either chain recurrent or “gradient-like”. 

\begin{theorem} \label{CONLEY}
$(i)$ Let $f: X \to X$ be a homeomorphism on a compact metric space. There exists a continuous Lyapunov function $u: X \to  \mathbb{R}$ for $f$ such that $N(u) = \mathcal{CR}(f)$. \\
$(ii)$ Let $\phi_t: X \times \mathbb{R} \to X$ be a continuous flow on a compact metric space. There exists a continuous Lyapunov function $u: X \to  \mathbb{R}$ for $\phi_t$ such that $N(u) = \mathcal{CR}(\phi_t)$. \end{theorem}

\begin{figure}[h!]
\centering
  \centering
\includegraphics[width=.6\linewidth]{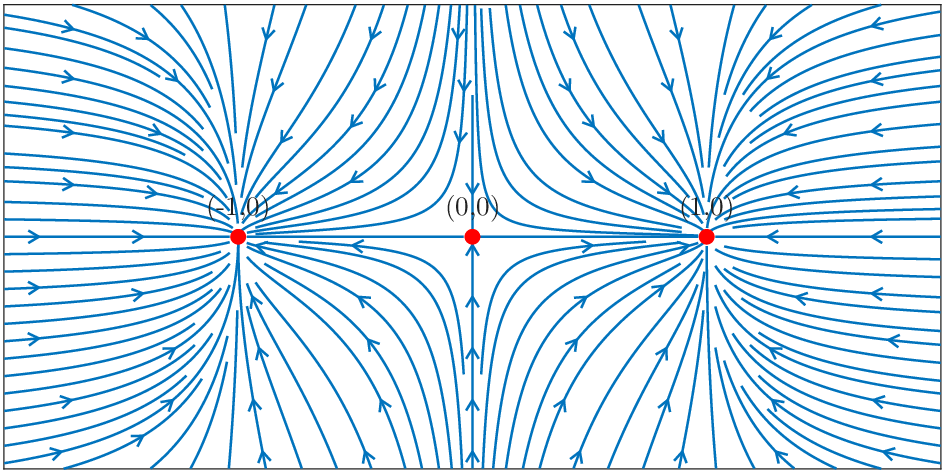}
\caption{The flow of Example \ref{example2.5}, together with its equilibria $(-1,0)$, $(0,0)$ and $(1,0)$.
\label{fig:ex_2.5}}\end{figure}

\begin{figure}[h!]
\centering
\begin{subfigure}{.15\textwidth}
  \centering
   \begin{tikzpicture}
  \node at (0,0) { \includegraphics[width=\linewidth]{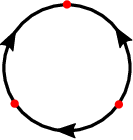}};
\node at (-0.7,-0.5) {$\textcolor{red}{a}$};
\node at (0.7,-0.5) {$\textcolor{red}{c}$};
\node at (0,0.9) {$\textcolor{red}{b}$};
  \end{tikzpicture}  
  \caption{}
  \label{fig:ex26a}
\end{subfigure}\hspace{1.5cm}
\begin{subfigure}{.15\textwidth}
  \centering
     \begin{tikzpicture}
  \node at (0,0) { \includegraphics[width=\linewidth]{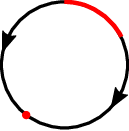}};
\node at (-0.6,-0.6)  {$\textcolor{red}{a}$};
\node at (0.8,0.5) {$\textcolor{red}{c}$};
\node at (0,0.9) {$\textcolor{red}{b}$};
  \end{tikzpicture}  
  \caption{}
  \label{fig:ex26b}
\end{subfigure}\\
\caption{The flow $\phi_t$ of Examples \ref{ex:2.6}(i) and \ref{ex:2.6}(ii) ((a) and (b), respectively).\label{fig:ex26}}
\end{figure}

\noindent Conley showed that in fact $u$ is a \emph{complete Lyapunov function}, that is -- in addition to decomposing the space into a chain recurrent part and a “gradient-like” part -- it possesses the following two properties: 
\begin{enumerate}
\item[$(a)$] If $x,y \in \mathcal{CR}$, then $u(x) = u(y)$ if and only if $x$ and $y$ belong to the same ``chain transitive component'' of $\mathcal{CR}$;
\item[$(b)$] $u(\mathcal{CR})$ is a compact nowhere dense subset of $\mathbb{R}$.
\end{enumerate}
\noindent In what follows, we offer just a glimpse of the proof of Theorem \ref{CONLEY} by highlighting a few key details: the decomposition of the chain recurrent set in terms of attractors and a brief look at the explicit construction of the complete Lyapunov function, providing the decomposition of the space. This construction gives the Lyapunov function with neutral set $\mathcal{CR}$; the additional completeness properties follow from the refined construction in Conley’s theorem. We give this brief overview for a flow, however the same underlying ideas provide the proof in the case of a homeomorphism. \\
~\newline 
\noindent Studying the qualitative properties of dynamical systems inevitably leads to consideration of objects called attractors and repellers. Here and in the rest of the manuscript, we use ${\rm{cl}}(Y)$ to indicate the closure of a given set $Y$.
\begin{definition} Let $\phi_t: X \times \mathbb{R} \to X$ be a continuous flow. A subset $A \subset X$ is called an \emph{attractor} if it is the $\omega$-limit of one of its neighborhoods. That is, if there exists a neighborhood $U$ of $A$ such that
$$\omega(U) = \bigcap_{T \ge 0} {\rm{cl}}(\phi_{[T,+\infty)}(U)) = A .$$
\end{definition}
\noindent The $\omega$-limit is easily seen to be closed and invariant. For compact ambient spaces, $\omega(U)$ is in particular compact and non-empty. If in addition $U$ is connected, then $\omega(U)$ is connected; see \cite[p. 88]{Ding2005} for detailed proofs. We need to remark that, in general: 
$$\bigcup_{x \in U} \omega(x) \subsetneq \omega(U) ,$$
as explained  below.
\begin{example}\label{example2.5} 
Consider the vector field on (the non compact space) $\mathbb{R}^2$:
$$\begin{cases}
\dot{x} = x(1-x^2), \\
\dot{y} = -y.
\end{cases}
$$
See Figure \ref{fig:ex_2.5} for a visualization of the dynamics. If $U := (-2,2) \times (-1,1)$, then it holds that $\omega(U) = [-1,1] \times \{0\}$, whereas $\bigcup_{(x,y) \in U} \omega((x,y)) = \{(\pm 1,0),(0,0)\}$.
\end{example}
\noindent In the case where $(X,d)$ is a compact metric space, associated with any attractor $A \subsetneq X$ there is a non-empty set called the \emph{complementary repeller} defined as follows: 
$$A^* = \{x \in X: \ \omega(x) \cap A = \emptyset\} ,$$
which is an attractor for the reverse flow. Moreover,
$$C(A,A^*) := X \setminus (A \cup A^*)$$
is called the set of \emph{connecting orbits} of the pair $(A,A^*)$. These definitions will be clarified in the next example.
\begin{example} \label{ex:2.6}
$(i)$ On the circle $\mathbb{S} = \mathbb{R}/\mathbb{Z}$ equipped with the standard quotient metric, consider the flow $\phi_t$ of Figure \ref{fig:ex26a}, where the bold red line and the arrows denote respectively fixed points and the direction of $\phi_t$. This dynamical system has two attractors. The first one is $A = \{b\}$, whose
corresponding $A^* = {\rm{cl}}(ca)$ (where $ca$ is the arc in $\mathbb{S}$ clockwise joining $c$ and $a$) and $C(A,A^*) = ab \cup bc$. The other attractor is $B = {\rm{cl}}(ab)$, whose corresponding $B^* = \{c\}$ and $C(B,B^*) = ca \cup bc$.
\\
$(ii)$ The dynamical system of Figure \ref{fig:ex26b} has only one attractor $A = \{a\}$, with $A^* = {\rm{cl}}(bc)$ and $C(A,A^*) = ab \cup ca$.
\end{example}

\noindent The first part of Conley's proof rests on decomposing the recurrent set using attractor-repeller pairs. In particular, it holds that $$\mathcal{CR}(\phi_t) = \bigcap\{ A \cup A^*: \, A \text{ is an attractor}\},$$ see \cite[p. 37]{Conley1978}. Moreover, as a consequence of the compactness of $X$, the set of attractors is at most countable, so that
\begin{equation} \label{first step}
\mathcal{CR}(\phi_t) = \bigcap\{ A_n \cup A_n^*: \, A_n \text{ is an attractor}\}.
\end{equation}
The second part of Conley's proof is based on the explicit construction of a Lyapunov function strictly decreasing on the connecting orbit of an attractor-repeller pair. In particular, he proves that for any attractor-repeller pair $(A, A^*)$, there exists a continuous function $h : X \to \mathbb{R}$ such that:
\begin{enumerate}
\item[$(a)$] $h^{-1}(0) = A$; 
\item[$(b)$] $h^{-1}(1) = A^*$;
\item[$(c)$] $h$ is strictly decreasing on the connecting orbit $C(A,A^*)$.
\end{enumerate}
The construction of $h$ follows these main lines. Define first 
$$\ell(x) := \frac{d(x,A)}{d(x,A) + d(x,A^*)}$$
and then 
$$k(x) := \sup_{t \ge 0} \{ \ell(\phi_t(x)) \}.$$
Since $A$ and $A^*$ are invariant, it holds that $k^{-1}(0) = A$, $k^{-1}(1) = A^*$, $k(X) \subseteq [0,1]$ and $k$ is decreasing on $C(A,A^*)$. The slightly technical part is proving that $k$ is continuous. The desired continuous function $h : X \to \mathbb{R}$, strictly decreasing on $C(A,A^*)$, is then given by
\begin{equation} \label{second step}
h(x) := \int^{+\infty}_0 e^{-s} k(\phi_s(x))ds .
\end{equation}
With these two steps, the Fundamental Theorem of Dynamical Systems is almost established. In fact, let $h_1, h_2, \ldots$ be the Lyapunov functions for the attractor-repeller pairs $(A_1,A_1^*), (A_2,A_2^*), \ldots$ of the flow $\phi_t$, as constructed in (\ref{second step}). We define $u: X \to \mathbb{R}$ as
$$u(x) := \sum_{n = 1}^{+\infty}\frac{h_n(x)}{3^n} .$$
Since every $h_n$ is continuous, $u$ is as well. Moreover, by the properties of the functions $h_n$ and the decomposition (\ref{first step}), $u$ is strictly decreasing along orbits outside the chain recurrent set and constant on orbits in the chain recurrent set. \\ 
~\newline
\noindent It is easy to construct a dynamical system whose chain recurrent set is the whole ambient space. In such a case, Conley's theorem gives a first integral. We propose below two examples where, even if $\mathcal{CR} = X$, part of the space is ``gradient like'', i.e. the underlying dynamics admits a Lyapunov function which is not a first integral. Here and in the rest of the manuscript, we use $\rm{Fix}$ and $\rm{Per}$ to indicate fixed points and periodic points, respectively. 

\begin{figure}[h!]
\centering\centering
\includegraphics[width=.6\linewidth]{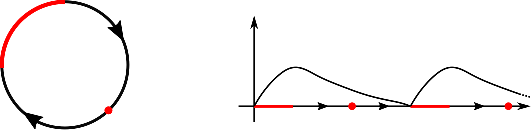}
\caption{{The flow $\phi_t$ of Example \ref{primo esempio}$(i)$ and the corresponding Lyapunov function strict outside ${\rm{Fix}}(\phi_t)$.
\label{fig:fig2}}}
\end{figure}

\begin{figure}[h!]
\centering
  \centering
\begin{tikzpicture}
 \node at (0,0) {\includegraphics[width=.35\linewidth]{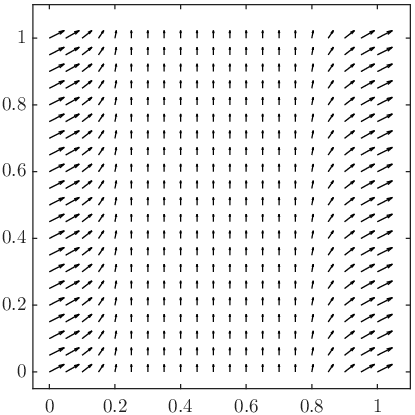}};
  \end{tikzpicture}
\caption{The flow of Example \ref{primo esempio}$(ii)$.
\label{fig:flow_torus}}\end{figure}

\begin{example} \label{primo esempio}
$(i)$ On $\mathbb{S} = \mathbb{R}/\mathbb{Z}$, let consider the flow $\phi_t$ of Figure \ref{fig:fig2}, where the bold line and the arrows denote respectively fixed points and the direction of $\phi_t$. In such a case, $\mathcal{CR}(\phi_t) = \mathbb{S}$ but $\phi_t$ clearly admits a Lyapunov function strictly decreasing outside ${\rm{Fix}}(\phi_t)$. \\
$(ii)$ On the torus $\mathbb{T}^2 = \mathbb{R}^2/\mathbb{Z}^2$ equipped with the standard quotient metric, consider the flow $\phi_t$ associated to the vector field
$$V(x,y) = 
(f(x),1)\, , 
$$
where \(f\) is a smooth nonnegative periodic function whose zero set is exactly \([1/4,3/4]\); see Figure \ref{fig:flow_torus}. In such a case, $\mathcal{CR}(\phi_t) = \mathbb{T}^2$ but $\phi_t$ clearly admits a Lyapunov function strictly decreasing outside ${\rm{Per}}(\phi_t) = [1/4,3/4] \times \mathbb{S}$.
\end{example}

\noindent The previous examples address the question of the existence of a set, strictly contained in $\mathcal{CR}$, whose complement is ``gradient-like''. This set, called \emph{strong chain recurrent set}, was introduced by R. Easton in 1978 (the same year of Conley's Fundamental Theorem) by replacing the requirement that \emph{each jump} be arbitrarily small -- see Definition~\ref{CR DEF}$(i)$ and $(ii)$ -- with the requirement that \emph{the sum of the jumps} be arbitrarily small. We will introduce and discuss strong chain recurrence in the next section. 

\section{Strengthening chain recurrence}\label{sec:strength}
The notion below -- to the best of the authors' knowledge -- first appeared in \cite{easton1978chain}. 
\begin{definition} \label{SCR epsilon}
$(i)$ Let $f: X \to X$ be a homeomorphism. A point $x \in X$ is said to be \emph{strong chain recurrent} for $f$ if for any $\varepsilon > 0$ there exists a finite sequence $(x_i)_{i=1,\dots,n} \subset X$ such that $x_1 =x$ and setting $x_{n+1} :=x$, we have 
\begin{equation} \label{omeo1}
\sum_{i = 1}^n d(f(x_i),x_{i+1}) < \varepsilon . 
\end{equation}
$(ii)$ Let $\phi_t: X \times \mathbb{R} \to X$ be a continuous flow. A point $x \in X$ is said to be \emph{strong chain recurrent} for $\phi_t$ if for any $\varepsilon > 0$ and $T > 0$, there exists a finite sequence $(x_i,t_i)_{i=1,\dots,n} \subset X \times \mathbb{R}$ such that $t_i \ge T$ for all $i$, $x_1 =x$ and setting $x_{n+1} :=x$, we have 
\begin{equation} \label{flusso1}
\sum_{i = 1}^n d(\phi_{t_i}(x_i),x_{i+1}) < \varepsilon .
\end{equation}
\end{definition}

\noindent The set of strong chain recurrent points is denoted by $\mathcal{SCR}$. For a continuous flow on a compact metric space, the strong chain recurrent set is invariant and closed; see \cite[p.~1265, after Def.~2.5]{bernardi2019conley}. Moreover, it is non-empty. In general, $\mathcal{SCR} \subseteq \mathcal{CR}$, and the inclusion can be strict. Finally, as clarified in the next example, $\mathcal{SCR}$ depends on the metric in $(X,d)$. For this reason, from now on, we use the notation $\mathcal{SCR}_d$. The next example is 
from \cite{Fathi_Pageault_2015_AubryMather}.

\begin{example} \label{metrica}
On $\mathbb{S}$ with the usual topology, let consider a flow which fixes a Cantor set and moves all other points clockwise. Choose a metric $d$ for which the Cantor set has Lebesgue measure $0$. Call the resulting flow $\phi_t: \mathbb{R} \times (\mathbb{S},d) \to (\mathbb{S},d)$. Then $\mathcal{SCR}_{d}(\phi_t) = \mathbb{S}$. Alternatively, choose a metric $\tilde{d}$ for which the Cantor set has positive Lebesgue measure. Call the resulting flow $\tilde{\phi}_t: \mathbb{R} \times (\mathbb{S},\tilde{d}) \to (\mathbb{S},\tilde{d})$. Then $\mathcal{SCR}_{\tilde{d}}(\tilde{\phi}_t) = {\rm{Fix}}(\tilde{\phi}_t)$.
\end{example}
\noindent The fixed points of Example \ref{primo esempio}$(i)$ and the periodic points of Example \ref{primo esempio}$(ii)$ are exactly the strong chain recurrent set for the corresponding flows. This means that in such examples there exists a Lyapunov function whose neutral set coincides with the strong chain recurrent set. It is therefore natural to ask whether a Conley-type theorem also holds for $\mathcal{SCR}_d$. The first affirmative answer to this question dates back to several years after the ``Fundamental Theorem of Dynamical Systems''. In fact, it is only in 2015 that A. Fathi and P. Pageault proved the following theorem in the discrete case.  
\begin{theorem} \label{FP}
Let $f: X \to X$ be a homeomorphism on a compact metric space $(X,d)$. There exists a Lipschitz-continuous Lyapunov function $u: X \to  \mathbb{R}$ for $f$ such that $N(u) = \mathcal{SCR}_d(f)$.
\end{theorem}
\noindent This result -- see \cite[Prop. 4.1]{pageault2011lyapunov} and \cite[Prop. 2.3]{Fathi_Pageault_2015_AubryMather} -- clearly corresponds to point $(i)$ of Theorem \ref{CONLEY}, with $\mathcal{SCR}_d$ in place of $\mathcal{CR}$. We emphasize that there is a gain in regularity, since in Theorem \ref{FP} the Lyapunov function is Lipschitz-continuous by construction. As in the case of Conley's original theorem, we believe it is didactically useful to outline the main ideas of the proof of Theorem \ref{FP}. The interesting fact is that Fathi and Pageault use techniques inspired by Weak KAM Theory, known for its applications in areas far from Topological Dynamics; especially in Aubry--Mather Theory and in the study of the Hamilton--Jacobi equation. \\
~\newline
\noindent First of all, let us rephrase the problem from a different point of view. Given $x,y \in X$, let $C(x,y)$ be the set of chains from $x$ to $y$, that is
$$C(x,y) = \{(x_i)_{i=1,\ldots,n+1}: x_1 = x, x_{n+1} = y\} .$$ 
Each element $C \in C(x,y)$ provides a ``strategy'' to go from $x$ to $y$, to which we associate a ``cost" given by:
$$\ell_d(C) := \sum_{i = 1}^n d(f(x_i),x_{i+1}) .$$
We now observe that the parameter $\varepsilon > 0$ in Definition \ref{SCR epsilon}$(i)$ is arbitrarily small. This suggests that we define the function
\begin{align*}
L_d:X \times X \to& [0,+\infty) \\
(x,y) \mapsto& \inf\{\ell_d(C): C \in C(x,y)\} ,
\end{align*}
so that, within this formulation:
$$\mathcal{SCR}_d(f) = \{x \in X: L_d(x,x) = 0\} .$$
The function $L_d$ has the following properties:
\begin{enumerate}
\item $L_d(x, y) \le L_d(x, z) + L_d(z, y)$ for every $x,y,z \in X$;
\item $L_d(x, f(x)) = 0$ for every $x \in X$;
\item $L_d(x, x)=0 \Leftrightarrow L_d(f(x), x)=0 \Leftrightarrow L_d(f(x), f(x)) = 0$ for every $x \in X$;
\item $|L_d(x, y) - L_d(x, z)| \le d(y, z)$ and $|L_d(x, y) - L_d(z, y)| \le d(f(x), f(z))$ for every $x,y,z \in X$. In particular, $L_d$ is continuous.
\end{enumerate}
With the previous properties, for the proofs of which we refer to \cite[Prop. 2.1]{Fathi_Pageault_2015_AubryMather} or \cite[Prop. 3.1]{pageault2011lyapunov}, Theorem \ref{FP} is almost proved. In fact, the last step is the following proposition, whose proof we include for the reader's convenience.
\begin{proposition} Let $(x_n)_{n \in \mathbb{N}}$ be a dense sequence in $X$. Then,
$$u(x) = \sum_{n = 1}^{+\infty} \frac{L_d(x_n,x)}{3^n}$$
is a Lipschitz-continuous Lyapunov function for $f$ such that $N(u) = \mathcal{SCR}_d(f)$.
\end{proposition}
\noindent \textit{Proof.} The function $u$ is well-defined and continuous because $L_d$ is bounded on the compact set $X \times X$. It is actually Lipschitz-continuous, in fact. Using properties 1 and 4 above, we obtain:
$$
|u(y) - u(x)| = \left|\sum_{n=1}^{+\infty} \frac{1}{3^n} ( L_d(x_n,y) - L_d(x_n,x)) \right| \le \frac{1}{2} d(x,y) ,
$$
for all $x,y \in X$. \\
The fact that $u$ is a Lyapunov function for $f$ is a straightforward consequence of properties 1 and 2. Indeed,
$$
u(f(x)) - u(x) = \sum_{n=1}^{+\infty} \frac{1}{3^n} ( L_d(x_n,f(x)) - L_d(x_n,x) ) \le \frac{1}{2} L_d(x,f(x)) = 0 ,
$$
for all $x \in X$. \\
It remains to show that $N(u) = \mathcal{SCR}_d(f)$. On the one hand, if $x \in \mathcal{SCR}_d(f)$ and using properties 1 and 3, we obtain:
$$u(x) - u(f(x)) = \sum_{n=1}^{+\infty} \frac{1}{3^n} ( L_d(x_n,x) - L_d(x_n,f(x)) ) \le \frac{1}{2} L_d(f(x),x) = 0 .$$
Since we already know that $u$ is a Lyapunov function for $f$, from the above inequality we conclude that $\mathcal{SCR}_d(f) \subseteq N(u)$. On the other hand, let $x \in X \setminus \mathcal{SCR}_d(f)$. Then we have that $L_d(x,x) > 0$. Corresponding to this $x \in X$, we introduce the continuous function $\mathcal{L}_d(y) := L_d(y,f(x)) - L_d(y,x) - L_d(x,f(x))$. By construction, $\mathcal{L}_d(x) = -L_d(x,x) < 0$. By the density of the sequence $(x_n)_{n \in \mathbb{N}}$, there exists at least one $n \in \mathbb{N}$ such that $\mathcal{L}_d(x_n) < 0$, or equivalently $L_d(x_n,f(x)) - L_d(x_n,x) < L_d(x,f(x))$. Consequently:
$$
u(f(x)) - u(x) = \sum_{n=1}^{+\infty} \frac{1}{3^n} ( L_d(x_n,f(x)) - L_d(x_n,x)) < \frac{1}{2} L_d(x,f(x)) = 0 .
$$
This proves that $N(u) \subseteq \mathcal{SCR}_d(f)$. We have then obtained the equality $N(u) = \mathcal{SCR}_d(f)$, which concludes the proof.
\hfill $\Box$ \\
~\newline
We remark that K. Yokoi independently proved the existence of a continuous (not \textit{a priori} Lipschitz-continuous) Lyapunov function for $f: X \to X$, strictly decreasing outside $\mathcal{SCR}_d(f)$ (see \cite[Thm. 5.2]{KatsuyaYokoi2015}). \\
~\newline
At this point, it is natural to ask whether the same result can be rephrased for a flow rather than a homeomorphism. The validity of this result has been demonstrated by O. Bernardi, A. Florio, and J. Wiseman in two separate papers, see \cite[Thm. 4.1]{bernardi2019existence} and \cite[Thm. 1.1]{bernardi2021conley}. \\
We recall that a continuous flow $\phi_t: X \times \mathbb{R} \to X$ on a compact metric space $(X, d)$ is uniformly Lipschitz continuous on compact subsets of $[0,+\infty)$ if for any $T > 0$ there exists $M_T > 0$ such that for every $x,y \in X$:
\begin{equation} \label{uni Lip cont}
d(\phi_t(x),\phi_t(y)) \le M_T d(x,y) \qquad \forall t \in [0,T] .
\end{equation}
\begin{theorem} \label{BFW}
$(i)$ Let $\phi_t: X \times \mathbb{R} \to X$ be a continuous flow on a compact metric space $(X,d)$ which is uniformly Lipschitz-continuous on compact subsets of $[0, +\infty)$. There exists a Lipschitz-continuous Lyapunov function $u: X \to  \mathbb{R}$ for $\phi_t$ such that $N(u) = \mathcal{SCR}_d(\phi_t)$. \\
$(ii)$ Let $\phi_t: X \times \mathbb{R} \to X$ be a continuous flow on a compact metric space $(X,d)$. There exists a continuous Lyapunov function $u: X \to  \mathbb{R}$ for $\phi_t$ such that $N(u) = \mathcal{SCR}_d(\phi_t)$.
\end{theorem}
\noindent The proof of Theorem \ref{BFW}$(i)$, given in \cite{bernardi2019existence}, combines the variational approach established by Fathi and Pageault in \cite{Fathi_Pageault_2015_AubryMather} and some of Conley’s techniques presented in \cite{Conley1978}. We stress that it applies to continuous flows which are uniformly Lipschitz-continuous on compact subsets of $[0, +\infty)$. We remark that, with the techniques adopted in \cite{bernardi2019existence}, the hypothesis on uniformly Lipschitz regularity of the flow with respect to time is due to the fact that, in building the Lyapunov function for the flow, some ``regularizing'' process of an initial function (coming from the variational approach) is needed. In particular, some Lipschitz-like control over time is required. \\
Even if any flow of a Lipschitz continuous vector field satisfies the regularity hypothesis of the above theorem, there are examples of (dynamically interesting) flows that fail to fulfill such a hypothesis; see e.g. the next Example \ref{es Martin}.
\begin{example} \label{es Martin}
Let be $\tau:[0,1]\rightarrow \mathbb{R}$ defined as follows.
$$
\tau(x):=\begin{cases}
\sqrt{\frac{1}{2}-x}+\frac{\sqrt{2}-1}{\sqrt{2}}\quad x\in\left[0,\frac{1}{2}\right] \\ \\
\sqrt{x-\frac{1}{2}}+\frac{\sqrt{2}-1}{\sqrt{2}}\quad x\in\left[\frac{1}{2},1\right].
\end{cases}
$$
On $X := \{ (x,y)\in \mathbb{R}^2 :\ x\in[0,1], \ 0\leq y\leq \tau(x) \}$ endowed with the standard metric from $\mathbb{R}^2$, identify the graph of $\tau$ with $[0,1]$: $(x,\tau(x))\sim(x,0)$, and the left side with the right side: $(0,y) \sim (1,y)$. Finally, consider the flow $\phi_t$ presented in Figure \ref{funz-tetto}. \\
We first observe that $\phi_t$ is a continuous flow, not uniformly Lipschitz-continuous on compact subsets of
$[0, +\infty)$. Indeed, since the function $\tau$ is not Lipschitz-continuous at $x = 1/2$, for any time $T$ greater than
the maximum of $\tau$, we cannot find a uniform constant $M_T$ such that the condition of uniform Lipschitz-continuity (\ref{uni Lip cont}) is satisfied. However, the construction of a Lyapunov function strict outside $\mathcal{SCR}_d(\phi_t)$ is instructive since $\mathcal{SCR}_d(\phi_t) \ne X$, as explained below. The set of periodic points corresponds to the central vertical strip (of positive measure), so that ${\rm{Per}}(\phi_t) \ne
X$. Moreover, every point in $X \setminus {\rm{Per}}(\phi_t)$ cannot be strong chain recurrent: its dynamics goes
``from left to right'' and, since the vertical strip of periodic points has positive measure, the strip cannot be
crossed by a finite number of jumps whose sum is arbitrarily small. Therefore, ${\rm{Per}}(\phi_t) = \mathcal{SCR}_d(\phi_t) \ne \mathcal{CR}(\phi_t) = X$.
\begin{figure}[h!]
				 \centering
   \begin{tikzpicture}
  \node at (0,0) {\includegraphics[width=.7\textwidth]{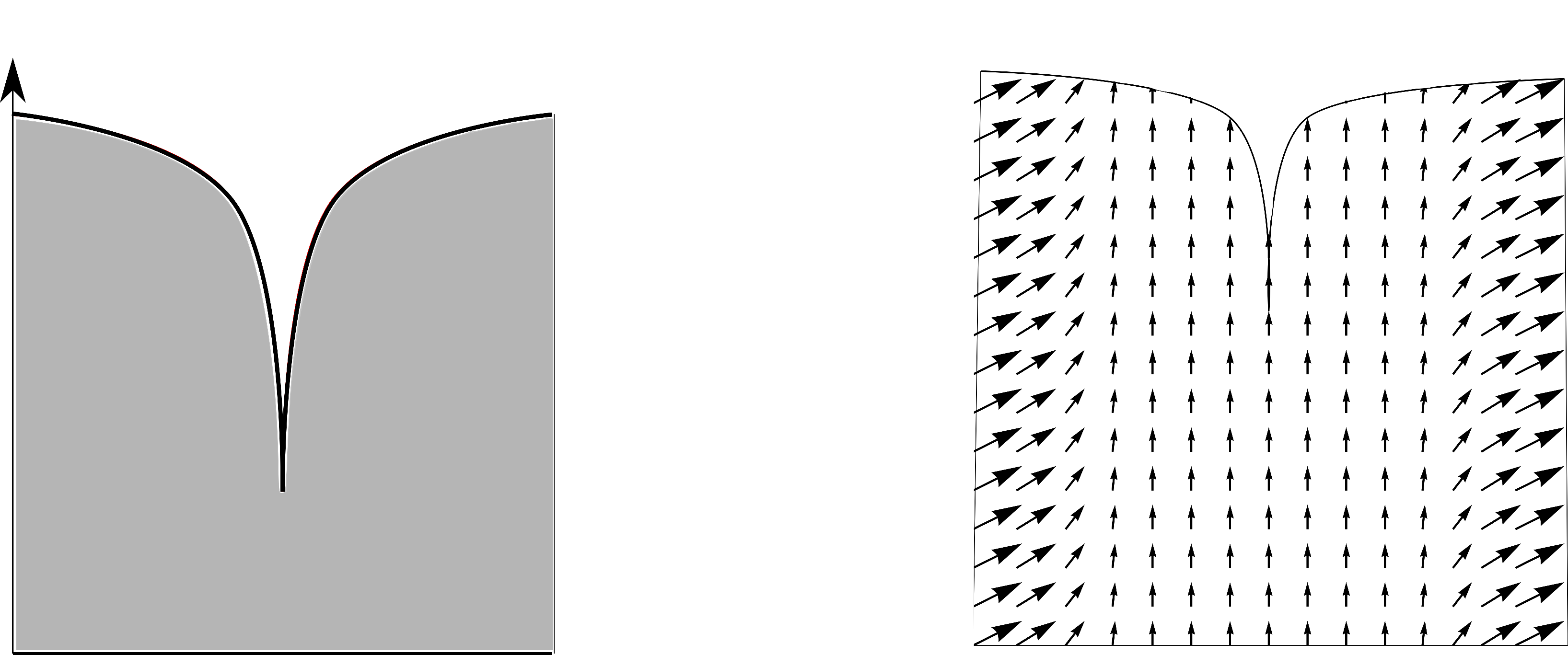}};
\node at (-5.6,-2.5) {$0$};
\node at (-1.5,-2.5) {$1$};
\node at (-1.6,1.75) {$\tau$};
  \end{tikzpicture}  
			\caption{The metric space and the flow $\phi_t$ of Example \ref{es Martin}.\label{funz-tetto}}
\end{figure}
\end{example}
\noindent The general case (Theorem \ref{BFW}$(ii)$), i.e.\ for merely continuous flows, is approached in \cite{bernardi2021conley} with arguments intimately close to Conley’s original ideas. Considering the purpose of the present work, we believe it is worth recalling the key points of the proof of Theorem \ref{BFW} as given in \cite{bernardi2021conley}, also because it involves a decomposition result of the strong chain recurrent set in the spirit of \eqref{first step}. To this aim, we need the definition of \emph{strongly stable set}.
\begin{definition} \label{s s set}
A closed set $B\subset X$ is \emph{strongly stable} if there exist a family $(U_{\eta})_{\eta\in(0,1)}$ of closed nested neighborhoods of $B$ and a function
	$$
	(0,1)\ni\eta\mapsto T(\eta)\in(0,+\infty)
	$$
	bounded on compact subsets of $(0,1)$, such that:
	\begin{itemize}
		\item[$(i)$] for any $0<\eta<\lambda<1$, $\{x\in X :\ d(x,U_{\eta})<\lambda-\eta\}\subseteq U_{\lambda}$; \vspace{2pt}
		\item[$(ii)$] $B=\bigcap_{\eta\in(0,1)}\omega(U_{\eta})$; \vspace{2pt}
		\item[$(iii)$] for any $0<\eta<1$, $\text{cl}\{ \phi_{[T(\eta),+\infty)}(U_{\eta}) \}\subseteq U_{\eta}$.
	\end{itemize}
\end{definition}
\noindent We refer to \cite[Def. 4.2]{bernardi2019conley}. Every strongly stable set $B$ is closed, forward invariant and stable, that is, it has a neighborhood basis of forward invariant sets; see \cite[Rmk. 4.1]{bernardi2019conley}. Moreover, every attractor is strongly stable; we refer to \cite[Ex. 4.1]{bernardi2019conley} for a detailed proof of this fact. Having established that an attractor is a strongly stable set, with the next Example \ref{SS not A} we prove that -- in general -- the converse is not true. In other words, Definition \ref{s s set} is really a generalization of the notion of an attractor.
\begin{example} \label{SS not A}
On the interval $[a,b]$, with the usual distance $d$, let us consider the dynamical system represented schematically in Figure \ref{fig:ex38}, where the bold line (and dots) and the arrows respectively denote fixed points and the direction of the flow $\phi_t$. In such a scenario, $\mathcal{SCR}_d(\phi_t)$ is the set of fixed points and the point $b$ is an attractor (and therefore strongly stable). However, every interval $[\alpha, \beta]$ of fixed points with $\beta<c$ is a strongly
stable set but not an attractor.
\end{example}
\begin{figure}[h!]
\centering
  \centering
\begin{tikzpicture}
 \node at (0,0) {\includegraphics[width=.35\linewidth]{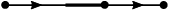}};
\node at (-2.7,0.4) {$a$};
\node at (2.7,0.4) {$b$};
\node at (0.7,0.4) {$c$};
  \end{tikzpicture}
\caption{Schematic representation of the dynamics described in Example \ref{SS not A}.
\label{fig:ex38}}\end{figure}

\noindent We are now ready to resume the main lines of the proof of Theorem \ref{BFW}. The first step involves the subtle relation between strongly stable sets and $\mathcal{SCR}_d(\phi)$: 
\begin{equation} \label{eccola}
	\mathcal{SCR}_d(\phi)=\bigcap \{ B\cup B^{\bullet} :\ B \text{ is strongly stable} \} ,
	\end{equation}
where 
$$B^\bullet := \{x \in X: \ \omega(x) \cap B = \emptyset\} .$$ 
We refer to \cite[Thm. 4.2]{bernardi2019conley} for the proof. \\
\noindent Secondly, for every pair $(B,B^\bullet)$ with $B$ strongly stable, it is possible to construct a continuous Lyapunov function $h$ for $\phi_t$ which is strictly decreasing on $X \setminus (B \cup B^\bullet)$. We note that, differently from Conley’s Lyapunov function for an attractor-repeller pair which is identically $0$ on the attractor and identically $1$ on the repeller, such a Lyapunov function $h$ assumes all the values between $0$ and $1$ in $B^\bullet$. \\
\noindent At this point, Conley’s proof of the Fundamental Theorem of Dynamical Systems uses the facts that attractor-repeller pairs are closed and that there are countably many such pairs. However, the pairs $(B,B^\bullet)$ involved in the decomposition of $\mathcal{SCR}_d(\phi_t)$, aren’t as well-behaved. Firstly, different $B$’s can give the same $B \cup B^\bullet$: the obvious example is the identity flow on a connected compact metric space. Moreover, $B \cup B^\bullet$ isn’t necessarily closed, there can be uncountably many different $B \cup B^\bullet$’s and, finally, the structure $B \cup B^\bullet$ is not closed under finite intersections. Although there are all these subtle differences with  Conley's case, the third step in the proof of Theorem \ref{BFW} is given by the existence of a -- at most countable -- collection of strongly stable sets $\{B_n\}$ such that
$$\mathcal{SCR}_d(\phi_t) = \bigcap B_n \cup B_n^{\bullet} .$$
We refer to \cite[Thm. 3.2]{bernardi2021conley}. \\
\noindent We finally conclude in the same way as explained in Section \ref{Charles Conley}. If $h_n$ is the Lyapunov function for the $n$th-pair $(B_n,B_n^\bullet)$, then
$$u(x) = \sum_{n = 1}^{+\infty}\frac{h_n(x)}{3^n}$$
is the desired continuous Lyapunov function for $\phi_t$, which is strictly decreasing outside $\mathcal{SCR}_d(\phi_t)$. 

~\newline
\noindent
We observe the following behavior of the strong chain recurrent set. 
For the flow $\phi_t$ in Example \ref{primo esempio}$(i)$, the strong chain recurrent set is the set of fixed points.  However, we can make an arbitrarily small perturbation to make the flow move the fixed points slightly clockwise, and then the strong chain recurrent set will be the entire circle $\mathbb S$.  Thus, there is an explosion in $\mathcal{SCR}$; that is, the mapping
\[
\psi_t \longmapsto \mathcal{SCR}(\psi_t)
\]
is not upper semicontinuous at $\psi_t=\phi_t$. Note further that there is no analogous explosion for the chain recurrent set. In fact, the chain recurrent set does not admit explosions; see, for instance, \cite{BlockFranke1985} and \cite[Thm.~2.2 and Rmk.~4.2]{bernardi2020generalized}. In contrast, explosions for the strong chain recurrent set can occur, as we see, and are discussed in \cite[Thm.~4.5]{bernardi2020generalized}. Implosions, on the other hand, may occur for both sets and are more common.

\section{Relations with neutral sets and generalized recurrence}\label{sec:neut_sets}
\begin{figure}[h!]
\centering
  \centering
\includegraphics[width=.85\linewidth]{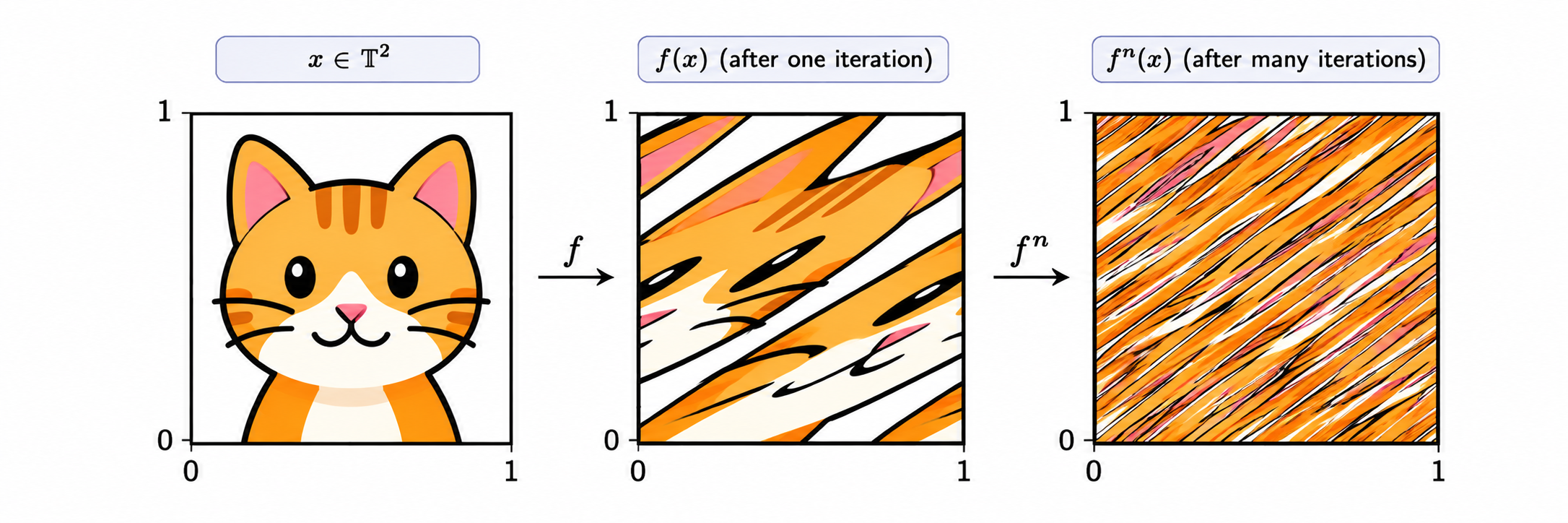}
\caption{Schematic illustration of Arnol’d’s cat map.}
\label{miao}\end{figure}

Let $f: X \to X$ be a homeomorphism on a compact metric space $(X,d)$ and 
$$\mathcal{L}_d(f) := \{\text{Lipschitz-continuous Lyapunov functions } u \text{ for } f \} .$$
\noindent Then, Theorem \ref{FP} clearly implies the following inclusion: 
$$\bigcap_{u \in \mathcal{L}_d(f)} N(u) \subseteq \mathcal{SCR}_d(f) .$$ 
Similarly, if $\phi_t: X \times \mathbb{R} \to X$ is a continuous flow on a compact metric space $(X,d)$ which is uniformly Lipschitz-continuous on compact subsets of $[0, +\infty)$, and 
$$\mathcal{L}_d(\phi_t) := \{\text{Lipschitz-continuous Lyapunov functions } u \text{ for } \phi_t \} ,$$
Theorem \ref{BFW}$(i)$ gives
$$\bigcap_{u \in \mathcal{L}_d(\phi_t)} N(u) \subseteq \mathcal{SCR}_d(\phi_t) .$$
It is important to note, however, that the previous inclusions are actually equalities, which consequently provide a further characterization of the strong chain recurrent sets in terms of neutral sets of the corresponding Lipschitz-continuous Lyapunov functions:
\begin{equation} \label{25 08 1}
\bigcap_{u \in \mathcal{L}_d(f)} N(u) = \mathcal{SCR}_d(f)   
\end{equation}
and
\begin{equation} \label{25 08 2}
\bigcap_{u \in \mathcal{L}_d(\phi_t)} N(u) = \mathcal{SCR}_d(\phi_t) .
\end{equation}
We refer respectively to \cite[§2.4]{pageault2011lyapunov}, \cite[Cor. 2.7]{Fathi_Pageault_2015_AubryMather}
and \cite[Thm. 2.2]{abbondandolo2018chain} for the proofs. \\
\noindent Taking into account the above characterization of the sets $\mathcal{SCR}_d$, it is therefore natural to consider the intersection of the neutral sets $N(u)$ as $u$ varies among all continuous (not necessarily Lipschitz-continuous) Lyapunov functions. In this way, we obtain the so-called \textit{generalized recurrent set} $\mathcal{GR}$. \\
~\newline
\noindent Generalized recurrence was originally introduced for flows in the 1960s by J. Auslander \cite{auslander1964generalized}, who defined the generalized recurrent set as the union of those orbits along which all continuous Lyapunov functions are constant. In more detail, we introduce the following definition.
\begin{definition}
Let $\phi_t: X \times \mathbb{R} \to X$ be a continuous flow  on a compact metric space $(X,d)$. The \emph{generalized recurrent set} $\mathcal{GR}(\phi_t)$ is defined as the intersection
$$\mathcal{GR}(\phi_t) := \bigcap_{u \in \mathcal{L}_0(\phi_t)} N(u) ,$$
where 
$$\mathcal{L}_0(\phi_t) := \{\text{continuous Lyapunov functions } u \text{ for } \phi_t \} .$$
\end{definition}
\noindent In the same paper, Auslander gave a characterization of this set in terms of the theory of prolongations. The generalized recurrent set was later extended to maps by Akin and Auslander (see \cite{akin1993general} and \cite{akin2010generalized}). More recently, Fathi and Pageault \cite{Fathi_Pageault_2015_AubryMather} proved that, for a homeomorphism, the generalized recurrent set can be equivalently defined by using Easton’s strong chain recurrence \cite{easton1978chain}, as stated in the following theorem.
\begin{theorem} Let $X$ be a compact metrizable topological space and $\mathcal{D}$ be the set of all metrics compatible with the topology of $X$. Let $f: X \to X$ be a homeomorphism. Then
$$\bigcap_{d \in \mathcal{D}} \mathcal{SCR}_d(f) = \mathcal{GR}(f) = \bigcap_{u \in \mathcal{L}_0(f)} N(u) ,$$
where
$$\mathcal{L}_0(f) = \{\text{continuous Lyapunov functions } u \text{ for } f \} .$$
\end{theorem}
\noindent The analogous result holds for flows, as shown in \cite[Prop.~6.4.1]{deleo2025streams}. 
We remark that -- with the theorem above -- we have truly introduced the smallest set on which every Lyapunov function is necessarily a first integral. On the one hand, we take into account the minimal regularity (continuity); on the other hand, we get rid of the metric (the basis of the definition of the strong chain recurrent) and we consider only the topology of $X$. \\
~\newline
\noindent We conclude this section with examples of generalized recurrent sets. 
\begin{example}
Consider the flow on $\mathbb S$ from Example~\ref{metrica}, which fixes a Cantor set and moves all other points clockwise.  Regardless of the choice of metric, the generalized recurrent set will be just the Cantor set of fixed points.  Similarly, the generalized recurrent set is just the fixed points in Example~\ref{primo esempio}$(i)$ and the periodic points in Example \ref{primo esempio}$(ii)$.
\end{example}
\begin{example}
Let $\phi_t$ be a flow on $\mathbb S$ which fixes a countable set of points and moves all other points clockwise.  Then $\mathcal{GR}(\phi) = \mathbb S$.
\end{example}
\begin{example}
Consider the doubling map on the circle $\mathbb S=\mathbb R/\mathbb Z$, $f(x) = 2x$, with periodic points 
\[
{\rm{Per}}(f) = \left\{\frac{k}{2^n-1}\pmod 1:
n\ge1,\ k\in\mathbb Z\right\}.
\] 
Every periodic point belongs to \(\mathcal{GR}(f)\); since \(\mathcal{GR}(f)\) is closed and periodic points are dense, \(\mathcal{GR}(f)=\mathbb S\). 
\end{example}
\begin{example}
Consider the so-called Arnol'd cat map on the torus $\mathbb T^2 = \mathbb R^2/\mathbb Z^2$:
$$f(x) = \left(\begin{matrix} 2 & 1 \\ 1 & 1 \end{matrix}\right) x\, .$$
The periodic points are dense, since they are the points with both coordinates rational, and so $\mathcal{GR}(f) = \mathbb T^2$. We refer to Figure \ref{miao} for a visualization.
\end{example}

\section{Relations with the stationary Hamilton--Jacobi equation}\label{sec:ham_jac}

Let $(M,g)$ be a closed, connected, smooth Riemannian manifold,
and let $H:T^*M\to\mathbb{R}$ be a $C^k$ Hamiltonian, with $k\geq 2$.
This section connects the previous results with the theory of
viscosity solutions of the stationary Hamilton--Jacobi equation
\begin{equation}\label{ViscoHJ}
H(x,du(x))=c ,
\end{equation}
where $c\in\mathbb{R}$ is a constant. Subsection~\ref{prima sotto} briefly reviews the notions of
viscosity solution, Ma\~n\'e critical value, and the projected Aubry set. Subsection~\ref{sotto seconda} introduces Ma\~n\'e Hamiltonians
and explores how the existence and regularity of Lyapunov functions
relate to those of viscosity subsolutions, with examples illustrating
the main ideas.

\subsection{Viscosity solutions, Ma\~n\'e critical value and the projected Aubry set} \label{prima sotto}
We start with the well-known notion of viscosity solution. 
\begin{definition}\label{def:sub_sup} Let $u \in C(M)$. \\
$(i)$ The function $u: M \to \mathbb{R}$ is a \emph{viscosity subsolution} of $H(x,d_x u) = c$ if, for every $\psi \in C^1(M;\mathbb{R})$ and every point $x_0 \in M$ such that $u-\psi$ has a local maximum at $x_0$, we have $H(x_0,d_{x_0} \psi) \le c$. \\
$(ii)$ The function $u: M \to \mathbb{R}$ is a \emph{viscosity supersolution} of $H(x,d_x u) = c$ if, for every $\psi \in C^1(M;\mathbb{R})$ and every point $x_0 \in M$ such that $u-\psi$ has a local minimum at $x_0$, we have $H(x_0,d_{x_0} \psi) \ge c$. \\
$(iii)$ The function $u:M \to \mathbb{R}$ is a \emph{viscosity solution} of $H(x,d_x u) = c$ if it is both a viscosity subsolution and a viscosity supersolution.
\end{definition}

\begin{figure}[h!]
\centering
  \centering
\includegraphics[width=.8\linewidth]{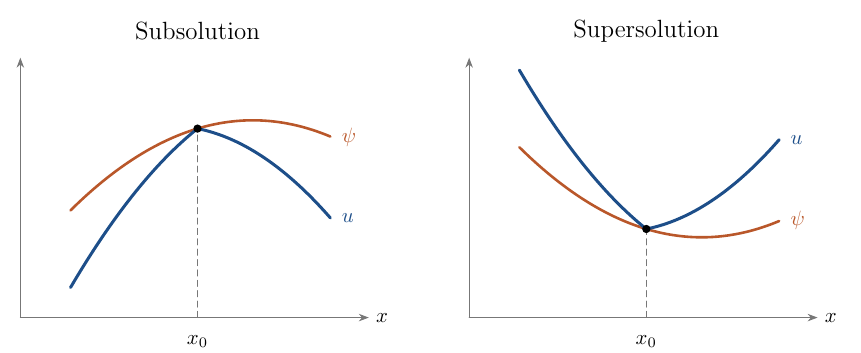}
\caption{Subsolutions and supersolutions.}
\label{sub-super-solution}
\end{figure}

\noindent We refer to Figure \ref{sub-super-solution} for a visualization of subsolutions and supersolutions. 

\noindent Definition \ref{def:sub_sup} gives meaning to the equation \eqref{ViscoHJ} at points where
$u$ is not necessarily differentiable, by testing it against smooth
functions that ``touch'' $u$ from above or below. It is worth noting that the terminology originates from the vanishing-viscosity method. Consider a family of classical solutions
$u^\varepsilon\in C^2(M)$ of
\begin{equation}\label{Viscmeth}
-\varepsilon\Delta_g u^\varepsilon
+H(x,du^\varepsilon(x))=c_\varepsilon ,
\end{equation}
where $\Delta_g=\operatorname{div}_g\nabla_g$.
If $c_\varepsilon\to c$ and $u^\varepsilon\to u$ uniformly
on $M$ as $\varepsilon\to0^+$, then $u$ is a viscosity
solution of $H(x,du(x))=c$.
This stability property does not require convergence of the
derivatives of $u^\varepsilon$. This is actually the motivation for this terminology in the original paper of M.G. Crandall and P.L. Lions \cite{crandall1983viscosity}.\\
~\newline
\noindent The study of viscosity solutions of the stationary Hamilton--Jacobi
equation associated with a Hamiltonian that is convex in the momentum
variable $p$, and more precisely satisfies the Tonelli assumptions\footnote{$D^2_{pp}H(x,p)$ positive definite, and
$H$ superlinear in $p$, uniformly with respect to $x$.},
naturally leads to two fundamental notions: the Ma\~n\'e critical
value $c(H)$ and the projected Aubry set associated with $H$.
\begin{definition}
The \emph{Ma\~n\'e critical value} of $H$ is
\begin{equation} \label{MANE}
c(H):=\inf_{u\in C^\infty(M)}\max_{x\in M}H(x,du(x)) .
\end{equation}
Equivalently, $c(H)$ is the smallest constant $c\in\mathbb R$
for which $H(x,du(x))\leq c$ admits a globally defined viscosity subsolution.
\end{definition}
\noindent Under these assumptions, $c(H)$ is the unique constant $c\in\mathbb R$
for which the stationary Hamilton--Jacobi equation
$H(x,du(x))=c$ admits a globally defined viscosity solution; see \cite[§~7 and~13]{Fathi2012}.
\\
\noindent To define the projected Aubry set, we introduce the Lagrangian
associated with $H$ by Legendre duality:
\[
L(x,v):=
\sup_{p\in T_x^*M}
\bigl\{p(v) - H(x,p)\bigr\} .
\]
For $t>0$ and $x,y\in M$, define the minimal action
\[
A_t(x,y):=
\inf_{\substack{\gamma\in\mathrm{AC}([0,t],M)\\
\gamma(0)=x,\ \gamma(t)=y}}
\int_0^t L\bigl(\gamma(s),\dot\gamma(s)\bigr)\,ds ,
\]
where $\mathrm{AC}([0,t],M)$ denotes the space of absolutely
continuous curves. The \emph{Peierls barrier} is
\[
h_H(x,y):=
\liminf_{t\to+\infty}\bigl(A_t(x,y)+c(H)t\bigr) .
\]

\begin{definition}
The \emph{projected Aubry set} associated with $H$ is
\[
\mathcal A(H):=
\bigl\{x\in M:h_H(x,x)=0\bigr\} .
\]
\end{definition}
\noindent For a Tonelli Hamiltonian on a compact manifold,
$\mathcal A(H)$ is a nonempty compact subset of $M$; see \cite[Thm. 5.2.8]{fathi_weak_kam}. 

\subsection{Lyapunov functions and viscosity critical subsolutions} \label{sotto seconda}

\begin{figure}[h!]
\centering
\begin{subfigure}{.15\textwidth}
  \centering
\includegraphics[width=\linewidth]{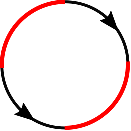}
  \caption{}
  \label{fig:ex55a}
\end{subfigure}\hspace{1.5cm}
\begin{subfigure}{.15\textwidth}
  \centering
\includegraphics[width=\linewidth]{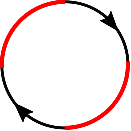}
  \caption{}
  \label{fig:ex55b}
\end{subfigure}\\
\caption{Dynamical system of Example \ref{41}\textit{(i)} and \ref{41}\textit{(ii)}.\label{figura2}}
\end{figure}

Let $V: M \to TM$ be a $C^k$ vector field, $k \ge 2$, inducing the flow $\phi_t: M \times \mathbb{R} \to M$. Denote by $\| v \|_x$ the norm of an element $v \in T_xM$ relative to the metric $g$ on $M$ and introduce the so-called Ma\~n\'e Hamiltonian
$$H_V: T^*M \to \mathbb{R}, \qquad H_V(x,p) = \frac{1}{2} \| p \|^2_x + p(V(x)) ,$$
where $\| p \|_x$ is the dual norm on $T_x^*M$ induced by $g$. We indicate by $\mathcal{A}_V$ the projected Aubry set associated to $H_V$. Since the constant functions are solutions of the Hamilton--Jacobi equation $H_V(x,du(x)) = 0$, the Ma\~n\'e critical level is $c(H_V) = 0$. \\
\noindent On the one hand, it holds that 
\begin{equation} \label{A in CR}
\mathcal{A}_V \subseteq \mathcal{CR}(\phi_t) ,
\end{equation}
see e.g. \cite[§I, Thm. V]{mane1997lagrangian}, \cite[Introduction, §3, Prop. 3.2]{pageault2011lyapunov} and \cite[Cor. 2]{Bernard2010}. \\
On the other hand, we know, by weak KAM theory, that there exists a $C^{1,1}$ critical subsolution $\tilde{u}:M \to \mathbb{R}$ of $H_V(x,du(x)) = 0$ which is strict outside $\mathcal{A}_V$; see \cite[Lemma 7]{Bernard2007}. In particular, since
$$H_V(x,d\tilde{u}(x)) \le 0 \Longrightarrow d\tilde{u}(V(x)) \le 0 ,$$ 
it follows that $\tilde{u}$ is a $C^{1,1}$ Lyapunov function for $\phi_t$, which is strict outside $\mathcal{A}_{V}$. This means that $N(\tilde{u})\subseteq\mathcal{A}_V$. \\
Recall also that the differentials of all $C^1$ critical
subsolutions agree on $\mathcal A_V$; see \cite[Lemma 6]{Bernard2007}. Since constant functions
are critical subsolutions, it follows that
\[
d\tilde u(x)=0
\qquad\text{for every }x\in\mathcal A_V.
\]
Moreover, $\mathcal A_V$ is invariant under $\phi_t$.
Consequently, for every $x\in\mathcal A_V$,
\[
\frac{d}{dt}\tilde u(\phi_t(x))
=
d\tilde u(\phi_t(x))(V(\phi_t(x)))
=0
\qquad\text{for every }t\in\mathbb R .
\]
Hence $\mathcal A_V\subseteq\mathcal N(\tilde u)$, and therefore
\[
\mathcal N(\tilde u)=\mathcal A_V.
\]
Since $\tilde u$ is Lipschitz-continuous on the compact manifold $M$,
equality~\eqref{25 08 2} implies
\begin{equation} \label{SCR in A}
\mathcal{SCR}(\phi_t)\subseteq
\mathcal N(\tilde u)=\mathcal A_V.
\end{equation}
\noindent From both inclusions \eqref{A in CR} and \eqref{SCR in A}, we conclude that in general
\begin{equation}\label{relazione Mane generale}
\mathcal{SCR}(\phi_t) \subseteq \mathcal{A}_V \subseteq \mathcal{CR}(\phi_t) ,
\end{equation}
so that the next proposition follows.  
\begin{proposition} \label{corollario mane}
If $\mathcal{SCR}(\phi)=\mathcal{CR}(\phi)$, then there exists a ${\cal{C}}^{1,1}$ Lyapunov function such that $
\mathcal{N}(\tilde{u}) =\mathcal{SCR}(\phi).
$  
\end{proposition}
\begin{proof}
If $\mathcal{SCR}(\phi_t) = \mathcal{CR}(\phi_t)$, then necessarily $\mathcal{SCR}(\phi_t) = \mathcal{A}_V = \mathcal{CR}(\phi_t)$ and the $C^{1,1}$ critical sub-solution $\tilde{u}:M \to \mathbb{R}$ of $H_V(x,d\tilde{u}(x)) = 0$ which is strict outside $\mathcal{A}_V$ gives the desired Lyapunov function with $N(\tilde{u}) =\mathcal{SCR}(\phi_t).$
\end{proof}
\noindent The next two examples discuss the result of Proposition \ref{corollario mane}. The basic Example \ref{41} shows two different smooth flows on the unit circle, having the following properties: 
\begin{itemize}
\item[$(i)$]  For the first flow, the hypothesis of Proposition \ref{corollario mane} holds and therefore there exists a differentiable Lyapunov function which is strict outside the strong chain recurrent set. 
\item[$(ii)$] For the second flow, we have that $\mathcal{SCR}(\phi_t) \subsetneq \mathcal{A}_V = \mathbb{S}$. In such a case, every critical subsolution given by weak KAM theory provides a first integral for $\phi_t$, but the dynamical system admits a differentiable Lyapunov function strict outside of $\mathcal{SCR}(\phi_t)$.
\end{itemize}
\begin{example} \label{41} $(i)$ On the circle $\mathbb{S} = \mathbb{R}/\mathbb{Z}$ equipped with the standard quotient metric, consider the dynamical systems of Figures \ref{fig:ex55a} and \ref{fig:ex55b}, where the bold line and the arrows denote the fixed points $\textnormal{Fix}(\phi_t)$ and the direction of the flow $\phi_t$, respectively. \\
Relative to the dynamical system of Figure \ref{fig:ex55a}, the hypothesis of Proposition \ref{corollario mane} is satisfied since $\mathcal{SCR}(\phi_t) = \mathcal{CR}(\phi_t) = \textnormal{Fix}(\phi_t)$ and there clearly exists a $C^{1,1}$ (also $C^{\infty}$) Lyapunov function whose neutral set is the strong chain recurrent set. \\
$(ii)$ Otherwise, in the case of the dynamical system of Figure \ref{fig:ex55b}, $\mathcal{SCR}(\phi_t) = \textnormal{Fix}(\phi_t)$ and $\mathcal{CR}(\phi_t) = \mathbb{S}$. Moreover -- see point $(i)$ of Theorem 4.8 and point $(i)$ of Lemma 4.14 in \cite{fathi2009hausdorff} --  the projected Aubry set associated to the corresponding Ma\~n\'e Hamiltonian coincides with $\mathcal{CR}(\phi_t) = \mathbb{S}$. Consequently, every critical subsolution given by weak KAM theory is a first integral but -- according to Theorem \ref{BFW}$(i)$ -- there clearly exists a Lipschitz continuous (even smooth) Lyapunov function strict outside $\textnormal{Fix}(\phi_t)$. We refer to \cite{akin1993general,FathiPageaultSmoothing} for general conditions under which a continuous Lyapunov function can be assumed to be smooth.
\end{example}

\noindent In the next example, we discuss another case in which critical subsolutions of weak KAM theory provide first integrals while the Lyapunov function of Theorem \ref{BFW}$(i)$ is not a first integral.
\begin{example}\label{es2} 
On $\mathbb{T}^2=\mathbb{R}^2/\mathbb{Z}^2$ endowed with the standard quotient metric, consider the dynamical system of Example \ref{primo esempio}$(ii)$; see Figure \ref{fig:flow_torus}. \\
\noindent In such a case, the strong chain recurrent set coincides with the set of periodic points, that is
$$\mathcal{SCR}(\phi_t) = \left[\frac{1}{4},\frac{3}{4}\right]\times\mathbb{S}$$
and $\mathcal{CR}(\phi_t) = \mathbb{T}^2$. Moreover -- by points $(i)$ of Theorem 4.8 and $(i)$ of Lemma 4.14 in \cite{fathi2009hausdorff} again -- the projected Aubry set associated to the Ma\~n\'e Hamiltonian $$\frac{1}{2} (p_1^2 + p_2^2) + f(x)p_1 + p_2$$ is $\mathcal{A}_V=\mathcal{CR}(\phi_t) = \mathbb{T}^2$. Consequently, every critical subsolution given by weak KAM theory is a first integral for $\phi_t$ while Theorem \ref{BFW}$(i)$ gives a Lyapunov function for $\phi_t$ which is strict outside $\left[\frac{1}{4},\frac{3}{4}\right]\times\mathbb{S}$.
\end{example}
\noindent We finally remark that the converse implication of Proposition \ref{corollario mane} does not in general hold true. For example, consider the dynamical system of Figure \ref{figura2}$(b)$. In such a case, it is easy to construct a smooth Lyapunov function which is strict outside $\mathcal{SCR}(\phi_t) = \textnormal{Fix}(\phi)$. However, $\mathcal{SCR}(\phi_t) \subsetneq \mathcal{CR}(\phi_t) = \mathbb{S}$. 

\section{Illustrative examples from mathematical biology models}\label{sec:examples}

\subsection{A singularly perturbed system of ODEs modeling epidemics}\label{sec:example1}

In order to illustrate the distinctions between these sets in a non-trivial case, we consider the reduced SIRWS system described in \cite{jardon2021geometric} (adapted from the model originally introduced in \cite{dafilis2014dynamical}). The population is partitioned into four compartments: $S$, susceptible individuals; $I$, infected and infectious individuals; $R$, recovered and immune individuals; and $W$, individuals with waning immunity, who may either completely lose immunity  and become susceptible again ($W\to S$) or receive an immunity boost after
contact with an infectious individual ($W\to R$).

\medskip

\begin{figure}[h!]
\centering
  \begin{tikzpicture}
  \node at (0,0) {\includegraphics[width=.5\linewidth]{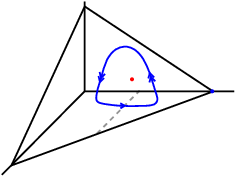}};
  \node at (0,1.6) {\textcolor{blue}{$\Gamma$}};
  \node at (3.4,0.2) {\textcolor{blue}{DFE}};
  \node at (0.6,0.6) {\textcolor{red}{EE}};
  \node at (4.2,0) {$S$};
  \node at (-4.2,-3) {$W$};
  \node at (-1.2,3.2) {$I$};
  \end{tikzpicture}
\caption{The regime $0<\epsilon\ll1$ for the three-dimensional system~\eqref{eq:SIRWS-3d}.\label{fig:ciclo_stab}}
\end{figure}

\noindent
The corresponding four-dimensional system is
\begin{equation}
\label{eq:SIRWS-4d}
\left\{
\begin{aligned}
\dot S &= -\beta SI+\epsilon\bigl(\xi(1-S)+2\kappa W\bigr),\\
\dot I &= \beta SI-\gamma I-\epsilon\xi I,\\
\dot R &= \gamma I+\nu\beta IW-\epsilon(2\kappa R+\xi R),\\
\dot W &= -\nu\beta IW+\epsilon(2\kappa R-2\kappa W-\xi W).
\end{aligned}
\right.
\end{equation}
Here $\beta$ is the infection rate, $\gamma$ is the recovery rate, $\epsilon\xi$ is the demographic
rate, $2\kappa$ is the rate governing the loss of waning immunity, $\nu\beta$ is the immunity
boosting rate nd $\epsilon$ represents the ratio between the timescales involved in the model.  Setting
\[
N:=S+I+R+W,
\]
we have
\[
\dot N=\epsilon\xi(1-N).
\]
Hence the hyperplane $N=1$ is invariant. Restricting to normalized initial conditions with
$N(0)=1$, we have $N(t)\equiv1$ and may eliminate
\[
R=1-S-I-W.
\]
We thus obtain the reduced three-dimensional system
\begin{equation}
\label{eq:SIRWS-3d}
\left\{
\begin{aligned}
\dot S &= -\beta SI+\epsilon\bigl(\xi(1-S)+2\kappa W\bigr),\\
\dot I &= \beta SI-\gamma I-\epsilon\xi I,\\
\dot W &= -\nu\beta IW
+\epsilon\bigl(2\kappa(1-S-I-2W)-\xi W\bigr).
\end{aligned}
\right.
\end{equation}
The biologically relevant phase space is the compact tetrahedron
\[
X:=\{(S,I,W)\in\mathbb{R}^{3}:S,I,W\geq0,\ S+I+W\leq1\}.
\]
\noindent Let $\phi_t$ denote the flow generated by the vector field in the ambient space. We restrict our analysis to complete orbits contained in the biologically relevant region $X$ and we compare the recurrent dynamics in two regimes: $0<\epsilon\ll1$ and the singular limit $\epsilon=0$.
~\newline
\noindent
\textbf{The case $0<\epsilon\ll1$.}
In \cite[§3.4.1, §3.4.2]{jardon2021geometric}, the authors show that, for suitable choices of the parameters, the system admits a stable limit cycle $\Gamma$.
In the dynamical regime represented in Figure~\ref{fig:ciclo_stab},
the disease-free equilibrium (DFE) is asymptotically stable for the dynamics restricted to
$\{I=0\}$, the endemic equilibrium (EE) is unstable, and the positive-infection dynamics is
attracted by $\Gamma$, apart from the invariant exceptional sets associated with the equilibria.
\noindent In this setting,
\begin{equation}
\label{eq:CR-SCR-positive-epsilon}
\mathcal{CR}(\phi_t)
=
\mathcal{SCR}(\phi_t)
=
\{\mathrm{DFE}\}\cup\{\mathrm{EE}\}\cup\Gamma.
\end{equation}
Indeed, since DFE and EE are equilibria and every point of $\Gamma$ is periodic,
\[
\{\mathrm{DFE}\}\cup\{\mathrm{EE}\}\cup\Gamma
\subseteq
\mathcal{SCR}(\phi_t)
\subseteq
\mathcal{CR}(\phi_t).
\]
Conversely, the asymptotic behavior described above excludes all remaining points from the
chain recurrent set.

\medskip

\begin{figure}[h!]
\centering
  \begin{tikzpicture}
  \node at (0,0) {\includegraphics[width=.5\linewidth]{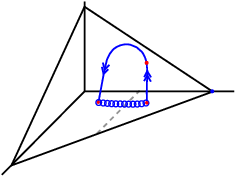}};
  \node at (0.4,1.7) {\textcolor{blue}{Heteroclinic Orbit}};
  \node at (3.4,0.2) {\textcolor{blue}{DFE}};
  \node at (1.2,0.9) {\textcolor{red}{A}};
  \node at (-0.95,-0.5) {\textcolor{red}{B}};
  \node at (1.3,-0.5) {\textcolor{red}{C}};
  \node at (4.2,0) {$S$};
  \node at (-4.2,-3) {$W$};
  \node at (-1.2,3.2) {$I$};
  \end{tikzpicture}  
\caption{The regime $\epsilon = 0$ for the three-dimensional system~\eqref{eq:SIRWS-3d} and an $(\varepsilon,T)$-chain from $A$ to itself.\label{fig:epsi0}}
\end{figure}

\noindent
\textbf{The singular case $\epsilon=0$.}
When $\epsilon=0$, system~\eqref{eq:SIRWS-3d} reduces to
\begin{equation}
\label{eq:SIRWS-singular}
\left\{
\begin{aligned}
\dot S &= -\beta SI,\\
\dot I &= I(\beta S-\gamma),\\
\dot W &= -\nu\beta IW.
\end{aligned}
\right.
\end{equation}
We assume $\beta>\gamma$ and set
\[
S^*:=\frac{\gamma}{\beta}\in(0,1). 
\] 
The line $S=S^*$ in the $SW$-plane is shown as a dashed line in Figures~\ref{fig:ciclo_stab} and \ref{fig:epsi0}.
The set (critical manifold, in Geometric Singular Perturbation Theory terminology) 
\[
\mathcal{C}_0:=\{(S,I,W)\in X:I=0\}
\]
consists entirely of equilibria. Moreover,
\[
\dot I=I(\beta S-\gamma),
\]
so that
\[
\mathcal{C}_0^r:=\{(S,0,W)\in \mathcal{C}_0:S>S^*\}
\]
is repelling in the $I$-direction, whereas
\[
\mathcal{C}_0^a:=\{(S,0,W)\in \mathcal{C}_0:S<S^*\}
\]
is attracting in the $I$-direction. \\
We now identify the non-constant complete orbits contained in $X$. Refer to Figure \ref{fig:epsi0} for a visualization of the following line of reasoning. Let $\mathcal H$ denote the union of all such orbits. Let $A=(S_0,I_0,W_0)\in X\setminus \mathcal{C}_0$ belong to a complete
orbit contained in $X$. Since $I_0>0$, uniqueness implies that $I(t)>0$ for every finite time
along this orbit. Moreover, necessarily $S(t)>0$ for every $t\in\mathbb R$. Indeed, if
$S=0$ and $I>0$, then
\[
\dot I=-\gamma I,
\]
and the corresponding backward orbit leaves the compact set $X$ in finite negative time. \\
Along the complete orbit,
\[
\dot S=-\beta SI<0,
\]
so that $S(t)$ is strictly decreasing. Two first integrals will be useful. First,
\[
\frac{d}{dt}
\left(
S+I-\frac{\gamma}{\beta}\log S
\right)=0,
\]
hence
\begin{equation}
\label{eq:first-integral-SI}
F(S,I):=S+I-S^*\log S
\end{equation}
is constant along the orbit. In particular, since $F$ is finite, $S(t)$ cannot tend to
zero. Secondly, a direct computation gives
\[
\frac{d}{dt}\left(\frac{W}{S^\nu}\right)=0,
\]
and therefore
\begin{equation}
\label{eq:first-integral-WS}
\frac{W}{S^\nu}=\mathrm{constant}
\end{equation}
along every complete orbit with $S>0$. \\
Since $S(t)$ is strictly decreasing and bounded, the limits
\[
S_-:=\lim_{t\to-\infty}S(t),
\qquad
S_+:=\lim_{t\to+\infty}S(t)
\]
exist and satisfy $S_->S_+>0$. Since the complete orbit is contained in the compact set $X$,
both its alpha- and omega-limit sets are nonempty. The strict monotonicity of $S$ implies
that every point in either limit set satisfies $I=0$. Moreover,
\eqref{eq:first-integral-SI} and \eqref{eq:first-integral-WS} imply that both limit sets
consist of a single point. Thus
\[
\alpha(A)=\{C\},
\qquad
\omega(A)=\{B\},
\]
where
\[
C=(S_-,0,W_-)\in \mathcal{C}_0,
\qquad
B=(S_+,0,W_+)\in \mathcal{C}_0.
\]
We now briefly show that
\[
S_->S^*>S_+.
\]
Indeed, evaluating the first integral~\eqref{eq:first-integral-SI} at the two endpoints gives
\[
S_- -S^*\log S_-
=
S_+ -S^*\log S_+.
\]
Let
\[
f(S):=S-S^*\log S.
\]
Then
\[
f'(S)=1-\frac{S^*}{S}
=\frac{S-S^*}{S}.
\]
Hence $f$ is strictly decreasing on $(0,S^*)$ and strictly increasing on $(S^*,+\infty)$,
and has its unique minimum at $S=S^*$. Since $S_->S_+$ and
$f(S_-)=f(S_+)$, the two values lie on opposite sides of $S^*$, namely
\[
S_->S^*>S_+.
\]
Consequently, every non-constant complete orbit contained in $X$ is a heteroclinic orbit
whose alpha-limit belongs to $\mathcal{C}_0^r$ and whose omega-limit belongs to $\mathcal{C}_0^a$.
Thus $\mathcal H$ is precisely the union of these heteroclinic orbits. \\
\noindent We now determine the chain recurrent set. Let
\[
M:=\bigcap_{t\geq0}\phi_t(X).
\]
The set $M$ is the maximal invariant subset of $X$, and the restriction of the dynamics to $M$ is a flow. Since $\mathcal{CR}(\phi_t)$ is itself invariant, every chain recurrent point belongs to $M$. By the discussion above,
\[
M=\mathcal{C}_0\cup\mathcal H.
\]
It follows that
\[
\mathcal{CR}(\phi_t)\subseteq \mathcal{C}_0\cup\mathcal H.
\]
The reverse inclusion follows from a constructive argument. Clearly, every point of $\mathcal{C}_0$ is
an equilibrium and therefore belongs to $\mathcal{CR}(\phi_t)$. Let now $A\in\mathcal H$.
Then, the complete orbit through $A$ is heteroclinic, with
\[
\lim_{t\to-\infty}\phi_t(A)=C\in \mathcal{C}_0^r,
\qquad
\lim_{t\to+\infty}\phi_t(A)=B\in \mathcal{C}_0^a.
\]
Fix $\varepsilon>0$ and $T>0$. Set $x_1:=A$ and $x_2:=B$. Choose
$t_1\geq T$ such that
\[
d\bigl(\phi_{t_1}(x_1),x_2\bigr)<\varepsilon.
\]
Since $\mathcal{C}_0$ is path connected, there exist finitely many points
\[
x_2=B,x_3,\ldots,x_{n-1}=C
\]
in $\mathcal{C}_0$ such that
\[
d(x_i,x_{i+1})<\varepsilon,
\qquad
i=2,\ldots,n-2.
\]
Since every $x_i\in \mathcal{C}_0$ is an equilibrium, for $i=2,\ldots,n-2$,
\[
d\bigl(\phi_T(x_i),x_{i+1}\bigr)
=
d(x_i,x_{i+1})
<
\varepsilon.
\]
Finally, since $C$ is the alpha-limit of $A$, we may choose a point $x_n$ on the same
heteroclinic orbit, such that
\[
d(C,x_n)<\varepsilon
\]
and such that the time $t_n$ needed to flow from $x_n$ to $A$ satisfies $t_n\geq T$.
Since $x_{n-1}=C$ is an equilibrium,
\[
d\bigl(\phi_T(x_{n-1}),x_n\bigr)
=
d(C,x_n)
<
\varepsilon,
\]
whereas, by construction,
\[
\phi_{t_n}(x_n)=A=x_1.
\]
Therefore, setting
\[
t_i:=T,
\qquad i=2,\ldots,n-1,
\]
and $x_{n+1}:=x_1=A$, we have
\[
d\bigl(\phi_{t_i}(x_i),x_{i+1}\bigr)<\varepsilon,
\qquad
i=1,\ldots,n.
\]
Thus $(x_i,t_i)_{i=1,\ldots,n}$
is an $(\varepsilon,T)$-chain from $A$ to itself. Since
$\varepsilon>0$ and $T>0$ are arbitrary, $
A\in\mathcal{CR}(\phi_t)$. Hence
\begin{equation}
\label{eq:CR-singular}
\mathcal{CR}(\phi_t)=\mathcal{C}_0\cup\mathcal H.
\end{equation}
\noindent The strong chain recurrent set is instead much smaller. Consider
\[
L(S,I,W):=S+I.
\]
Along the solutions of~\eqref{eq:SIRWS-singular},
\[
\frac{d}{dt}L(S,I,W)
=
\dot S+\dot I
=
-\gamma I.
\]
Thus $L$ is a smooth Lyapunov function for $\phi_t$, with strict decrease whenever $I>0$,
and therefore
\[
N(L)=\mathcal{C}_0.
\]
Since $L$ is smooth, hence Lipschitz-continuous on the compact set $X$, equality (\ref{25 08 2}) yields
\[
\mathcal{SCR}(\phi_t)
\subseteq
N(L)
=
\mathcal{C}_0.
\]
Conversely, every point of $\mathcal{C}_0$ is an equilibrium and hence strongly chain recurrent. Therefore
\begin{equation}
\label{eq:SCR-singular}
\mathcal{SCR}(\phi_t)=\mathcal{C}_0.
\end{equation}
\noindent Finally, $\mathcal H$ is nonempty. Indeed, the side $W=0$ is invariant. Fix
$S_c\in(S^*,1)$ and set
\[
C=(S_c,0,0)\in \mathcal{C}_0^r.
\]
On $W=0$, the level set
\[
F(S,I)=F(S_c,0)
\]
contains a non-constant branch with $I>0$ joining $C$ to a second point
\[
B=(S_b,0,0),
\qquad
S_b<S^*.
\]
Along this branch, $S+I$ is nonincreasing, so it remains in $X$. Hence this branch is a
heteroclinic orbit contained in $X$, and therefore belongs to $\mathcal H$. Combining \eqref{eq:CR-singular} and \eqref{eq:SCR-singular}, we obtain
\[
\mathcal{SCR}(\phi_t)
=
\mathcal{C}_0
\subsetneq
\mathcal{CR}(\phi_t)
=
\mathcal{C}_0\cup\mathcal H.
\]
In particular, this example gives a concrete mechanism producing a strict inclusion between
strong chain recurrence and chain recurrence. 

\subsection{A singularly perturbed system of ODEs for predator-prey models}\label{sec:example2}
Unlike the previous section, this example does not produce a distinction between \(\mathcal{CR}\) and \(\mathcal{SCR}\); rather, it illustrates how the common recurrent set changes through equilibria, a heteroclinic cycle, and a limit cycle as a single parameter of the model varies. In particular, we consider the two-dimensional system studied in \cite{borsotti2026geometric}. On
\[
\Delta=\{(u,v):0\le u\le1,\ v\ge0\},
\]
let
\begin{align} \label{eq2}
\left\{
\begin{aligned}
    \dot{u} &= \varepsilon u (u-l) (1-u)  - uv, \\ 
    \dot{v} &= -\gamma v(m-u), \\ 
\end{aligned}
\right.
\end{align}
where 
$$l\in (0,1), \quad \gamma>0, \quad m>0, \quad 0<\varepsilon \ll 1.$$ 
Let $\phi_t$ denote the flow generated by the vector field. We focus on the recurrent invariant sets contained in the biologically relevant region $\Delta$. This system models the interaction of resources $u$ and production $v$, and it can be interpreted as a predator-prey model, as well. The variables are normalized, so that $1$ is the carrying capacity of the resources, $l$ is the extinction threshold, $\gamma$ is the resource utilization efficiency and $\gamma m$ the decay rate of products. In particular, the growth of resources happens at a much slower rate than their consumption (represented by the small parameter $\varepsilon$). We refer to \cite{borsotti2026geometric} for a detailed explanation and derivation of the model.
\begin{figure}[h!]
\centering
\begin{subfigure}{.3\textwidth}
  \centering
  \begin{tikzpicture}
 \node at (0,0) {\includegraphics[width=\linewidth]{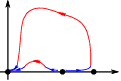}};
 \node at (-2.4,-1) {$\mathbf{x_1}$};
\node at (0.1,-1.6) {$\mathbf{x_2}$};
\node at (1.5,-1.6) {$\mathbf{x_3}$};
\node at (2.2,-1.6) {$u$};
\node at (-2.4,1.5) {$v$};
\node at (0,2) {Case 1: $m<l$};
  \end{tikzpicture}
  \caption{}
  \label{fig:sketch1}
\end{subfigure}\hspace{.5cm}
\begin{subfigure}{.3\textwidth}
  \centering
 \begin{tikzpicture}
 \node at (0,0) {\includegraphics[width=\linewidth]{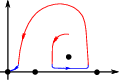}};
  \node at (-2.4,-1) {$\mathbf{x_1}$};
\node at (1.7,-1.6) {$\mathbf{x_3}$};
\node at (-1,-1.6) {$\mathbf{x_2}$};
\node at (0.5,-0.4) {$\mathbf{x_4}$};
\node at (2.2,-1.6) {$u$};
\node at (-2.4,1.5) {$v$};
\node at (0,2) {Case 2: $l<m<\tilde{m}(l,\gamma,\varepsilon)$};
  \end{tikzpicture}
  \caption{}
  \label{fig:sketch2}
\end{subfigure}\hspace{.5cm}
\begin{subfigure}{.3\textwidth}
  \centering
  \begin{tikzpicture}
 \node at (0,0) {\includegraphics[width=\linewidth]{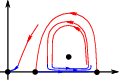}};
  \node at (-2.4,-1) {$\mathbf{x_1}$};
\node at (1.7,-1.6) {$\mathbf{x_3}$};
\node at (-1,-1.6) {$\mathbf{x_2}$};
\node at (0.5,-0.4) {$\mathbf{x_4}$};
\node at (2.2,-1.6) {$u$};
\node at (-2.4,1.5) {$v$};
\node at (0,2) {Case 3: $m=\tilde{m}(l,\gamma,\varepsilon)$};
  \end{tikzpicture}
  \caption{}
  \label{fig:sketch4}
\end{subfigure}\\
\begin{subfigure}{.3\textwidth}
  \centering
 \begin{tikzpicture}
 \node at (0,0) {\includegraphics[width=\linewidth]{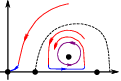}};
  \node at (-2.4,-1) {$\mathbf{x_1}$};
\node at (1.7,-1.6) {$\mathbf{x_3}$};
\node at (-1,-1.6) {$\mathbf{x_2}$};
\node at (0.5,-0.4) {$\mathbf{x_4}$};
\node at (2.2,-1.6) {$u$};
\node at (-2.4,1.5) {$v$};
\node at (0.75,0.9) {$v=\alpha(u)$};
\node at (0,2) {Case 4: $\tilde{m}(l,\gamma,\varepsilon)<m<\frac{l+1}{2}$};
  \end{tikzpicture}
  \caption{}
  \label{fig:sketch3}
\end{subfigure}\hspace{.5cm}
\begin{subfigure}{.3\textwidth}
  \centering
 \begin{tikzpicture}
 \node at (0,0) {\includegraphics[width=\linewidth]{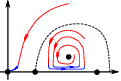}};
  \node at (-2.4,-1) {$\mathbf{x_1}$};
\node at (1.7,-1.6) {$\mathbf{x_3}$};
\node at (-1,-1.6) {$\mathbf{x_2}$};
\node at (0.4,-0.4) {$\mathbf{x_4}$};
\node at (2.2,-1.6) {$u$};
\node at (-2.4,1.5) {$v$};
\node at (0.75,0.9) {$v=\alpha(u)$};
\node at (0,2) {Case 5: $\frac{l+1}{2}\leq m<1$};
  \end{tikzpicture}
  \caption{}
  \label{fig:sketch5}
\end{subfigure}\hspace{.5cm}
\begin{subfigure}{.3\textwidth}
  \centering    \begin{tikzpicture}
 \node at (0,0) {\includegraphics[width=\linewidth]{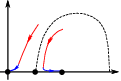}};
 \node at (-2.4,-1) {$\mathbf{x_1}$};
\node at (0.1,-1.6) {$\mathbf{x_3}$};
\node at (-1,-1.6) {$\mathbf{x_2}$};
\node at (2.2,-1.6) {$u$};
\node at (-2.4,1.5) {$v$};
\node at (1,1.4) {$v=\alpha(u)$};
\node at (0,2) {Case 6: $m>1$};
  \end{tikzpicture}
  \caption{}
  \label{fig:sketch6}
\end{subfigure}
\caption{Phase portrait for the system \eqref{eq2} as $m$ varies.}
\label{fig:sketch}
\end{figure}
\newline \noindent Consider Figure \ref{fig:sketch}, where the slow parts of the orbits are depicted in blue and with a single arrow along them, and the fast parts in red and with a double arrow. 
The dynamics illustrated therein can be summarized as follows. 

\medskip
\noindent 
\textbf{Case 1: $m<l$.} Every point in $\Delta$, except for $\mathbf{x}_2=(l,0)$ and points on the stable manifold of $\mathbf{x}_3=(1,0)$, is attracted to $\mathbf{x}_1=(0,0)$.

\medskip
\noindent
\textbf{Case 2: $l<m<\widetilde m(l,\gamma,\varepsilon)$.} The unstable equilibrium $\mathbf{x}_4=(m,\varepsilon(m-l)(1-m))$ appears. Every point in $\Delta$ outside the stable manifolds of $\mathbf{x}_2$ and $\mathbf{x}_3$, other than $\mathbf{x}_4$, is attracted to $\mathbf{x}_1$.

\medskip
\noindent
\textbf{Case 3: $m=\widetilde m(l,\gamma,\varepsilon)$.} The system exhibits a heteroclinic cycle connecting $\mathbf{x}_2$ to $\mathbf{x}_3$ along a slow orbit and $\mathbf{x}_3$ to $\mathbf{x}_2$ along a fast orbit. Every orbit starting strictly inside this cycle, except for the equilibrium $\mathbf{x}_4$, is attracted to it, whereas every orbit starting strictly outside the cycle is attracted to $\mathbf{x}_1$.

\medskip
\noindent
\textbf{Case 4:
$\widetilde m(l,\gamma,\varepsilon)<m<\frac{l+1}{2}$.} The stable manifold of $\mathbf{x}_2$ separates the basins of attraction of $\mathbf{x}_1$ and a unique stable limit cycle (purple) surrounding the unstable equilibrium $\mathbf{x}_4$. As $\varepsilon\to0$, this separatrix is approximated by the curve
\[
v=\alpha(u):=\gamma\left(l-u+m\log\left(\frac{u}{l}\right)\right),
\qquad l\leq u\leq1.
\]
Orbits in the interior of $\Delta$ on the limit-cycle side of this
separatrix, except for $\mathbf{x}_4$, are attracted to the cycle; those on the other side are attracted to $\mathbf{x}_1$. 

\medskip
\noindent
\textbf{Case 5.} For $(l+1)/2\leq m<1$, the stable manifold of $\mathbf{x}_2$ separates the basins of attraction of $\mathbf{x}_1$ and $\mathbf{x}_4$. Orbits starting in the interior of $\Delta$ below this separatrix are attracted to $\mathbf{x}_4$, whereas those starting above it are attracted to $\mathbf{x}_1$. The limit cycle shrinks to $\mathbf{x}_4$ as $m\to((l+1)/2)^{-}$.

\medskip
\noindent
\textbf{Case 6.} For $m>1$, every point in $\Delta$ outside the stable manifold of $\mathbf{x}_2$ is attracted either to $\mathbf{x}_1$ or to $\mathbf{x}_3$. This stable manifold separates the two basins of attraction and is again approximated by the curve $v=\alpha(u)$ as $\varepsilon\to0$.

In what follows, we discuss how $\mathcal{CR}(\phi_t)$ changes as $m$ increases. The bifurcation values \(m=l\) and \(m=1\) are not considered here. For the proof of these distinct behaviours, we once again refer to \cite{borsotti2026geometric}, focusing here on describing the results.

\medskip

\noindent 
\textbf{Case 1: $m<l$.}
In this regime, the relevant recurrent sets are only the equilibria
$\mathbf{x}_1,\mathbf{x}_2,\mathbf{x}_3$. All other trajectories exhibit the
asymptotic behavior described in Figure~\ref{fig:sketch}$(a)$, and there are no
periodic or heteroclinic cycles. Hence
\[
\mathcal{SCR}(\phi_t) = \mathcal{CR}(\phi_t)
=
\{\mathbf{x}_1,\mathbf{x}_2,\mathbf{x}_3\}.
\]

\medskip

\noindent
\textbf{Case 2: $l<m<\widetilde m(l,\gamma,\varepsilon)$.}
The interior equilibrium $\mathbf{x}_4$ lies in the relevant region, but there
are still no recurrent invariant sets different from equilibria. Therefore
\[
\mathcal{SCR}(\phi_t) = \mathcal{CR}(\phi_t)
=
\{\mathbf{x}_1,\mathbf{x}_2,\mathbf{x}_3,\mathbf{x}_4\}.
\]

\medskip

\noindent
\textbf{Case 3: $m=\widetilde m(l,\gamma,\varepsilon)$.}
At the critical value, the system exhibits a heteroclinic cycle connecting
$\mathbf{x}_2$ and $\mathbf{x}_3$. Let $\mathcal{H}_{\mathrm{het}}$
denote the whole heteroclinic cycle, including the equilibria
$\mathbf{x}_2$ and $\mathbf{x}_3$. Equivalently, if $H_{23}$ and $H_{32}$
denote the two open heteroclinic connections, then
\[
\mathcal{H}_{\mathrm{het}}
=
\{\mathbf{x}_2,\mathbf{x}_3\}\cup H_{23}\cup H_{32}.
\]
Every point of $\mathcal{H}_{\mathrm{het}}$ is chain recurrent, since one may
follow one heteroclinic branch close to the next equilibrium, make an
arbitrarily small jump to the other branch, and close the chain after following
the second branch, similarly to what was done in the case $\epsilon=0$ of Section \ref{sec:example1}. Moreover, the two jumps near $x_2$ and $x_3$ can be chosen so that the sum of their lengths is smaller than $\varepsilon$. Hence the same construction gives a strong chain, and every point of $\mathcal H_{\mathrm{het}}$ is strongly chain recurrent. Therefore,
\[
\mathcal{SCR}(\phi_t) = \mathcal{CR}(\phi_t)
=
\{\mathbf{x}_1,\mathbf{x}_4\}\cup\mathcal{H}_{\mathrm{het}},
\]
or, equivalently,
\[
\mathcal{SCR}(\phi_t) = \mathcal{CR}(\phi_t)
=
\{\mathbf{x}_1,\mathbf{x}_2,\mathbf{x}_3,\mathbf{x}_4\}
\cup H_{23}\cup H_{32}.
\]

\medskip

\noindent
\textbf{Case 4:
$\widetilde m(l,\gamma,\varepsilon)<m<\frac{l+1}{2}$.}
The heteroclinic cycle has disappeared and a stable periodic orbit,
denoted by $\Gamma_m$, is present around the unstable equilibrium
$\mathbf{x}_4$. Hence
\[
\mathcal{SCR}(\phi_t) = \mathcal{CR}(\phi_t)
=
\{\mathbf{x}_1,\mathbf{x}_2,\mathbf{x}_3,\mathbf{x}_4\}
\cup\Gamma_m.
\]

\medskip

\noindent
\textbf{Case 5:
$\frac{l+1}{2}\leq m<1$.}
The periodic orbit has collapsed onto the equilibrium $\mathbf{x}_4$, which is
now attracting. There are no additional recurrent non-equilibrium invariant
sets, and therefore
\[
\mathcal{SCR}(\phi_t) = \mathcal{CR}(\phi_t)
=
\{\mathbf{x}_1,\mathbf{x}_2,\mathbf{x}_3,\mathbf{x}_4\}.
\]

\medskip

\noindent
\textbf{Case 6: $m>1$.}
The equilibrium $\mathbf{x}_4$ no longer belongs to the relevant region, and
the only chain recurrent points are the equilibria on the $u$-axis:
\[
\mathcal{SCR}(\phi_t) = \mathcal{CR}(\phi_t)
=
\{\mathbf{x}_1,\mathbf{x}_2,\mathbf{x}_3\}.
\]
\noindent Summarizing,
\[
\mathcal{SCR}(\phi_t) = \mathcal{CR}(\phi_t)=
\begin{cases}
\{\mathbf{x}_1,\mathbf{x}_2,\mathbf{x}_3\},
& m<l,\\[2mm]
\{\mathbf{x}_1,\mathbf{x}_2,\mathbf{x}_3,\mathbf{x}_4\},
& l<m<\widetilde m(l,\gamma,\varepsilon),\\[2mm]
\{\mathbf{x}_1,\mathbf{x}_4\}\cup\mathcal{H}_{\mathrm{het}},
& m=\widetilde m(l,\gamma,\varepsilon),\\[2mm]
\{\mathbf{x}_1,\mathbf{x}_2,\mathbf{x}_3,\mathbf{x}_4\}\cup\Gamma_m,
& \widetilde m(l,\gamma,\varepsilon)<m<\dfrac{l+1}{2},\\[3mm]
\{\mathbf{x}_1,\mathbf{x}_2,\mathbf{x}_3,\mathbf{x}_4\},
& \dfrac{l+1}{2}\leq m<1,\\[3mm]
\{\mathbf{x}_1,\mathbf{x}_2,\mathbf{x}_3\},
& m>1.
\end{cases}
\]

\printbibliography

\end{document}